\documentstyle[12pt, twoside]{article}
\input amssym.def
\input amssym.tex

\newenvironment{demo}[1]%
{\vskip-\lastskip\medskip
  \noindent
  {\em #1.}\enspace
  }%
{\qed\par\medskip
  }

\newcommand{\qed}{
  \strut\hfill
  \mbox{$\Box$}
  }

\newcommand{\bn}{ \widetilde{H}_n }
\newcommand{\bG}{{H}\G}
\newcommand{\bgn}{ {\widetilde{H}\G}_n }
\newcommand{\C}{ {\Bbb C} }

\newcommand{\FG}{ {\cal F}_{\G}(X) }
\newcommand{\g}{\gamma}
\newcommand{\GG}{\widetilde{G}}
\newcommand{\G}{\Gamma}

\newcommand{\Gn}{\G_n}

\newcommand{\Hy}{\widetilde{H}}
\newcommand{\Hom}{ \mbox{Hom} }
\newcommand{\Ind}{\mbox{Ind} }

\newcommand{\kbgnfree}{ \underline{K}^-_{\bgn}(X^n) }
\newcommand{\la}{\lambda}

\newcommand{\NP}{\Pi}
\newcommand{\rk}{ {\rm rank }}
\newcommand{\sss}{ {Q}}
\newcommand{\tFG}{ {\cal F}^-_{\G}(X) }
\newcommand{\tFGvar}{ {\cal F}^-_{\G}(X,t) }
\newcommand{\tG}{ {\widetilde{H}\G}}

\newcommand{\tRG}{R^-_{\G}}
\newcommand{\Z}{ {\Bbb Z} }

\newtheorem{theorem}{Theorem}[section]
\newtheorem{lemma}{Lemma}[section]
\newtheorem{remark}{Remark}[section]

\newtheorem{definition}{Definition}[section]

\newtheorem{proposition}{Proposition}[section]
\newtheorem{corollary}{Corollary}[section]

\begin{document}
\title{Equivariant $K$-theory, generalized symmetric products, and twisted Heisenberg algebra}
\author{Weiqiang Wang\thanks{Partially supported by an NSF grant and an
FR\&PD grant at NCSU.}
}
\date{}
\maketitle

\begin{abstract}{
For a space $X$ acted by a finite group $\G$, the product space
$X^n$ affords a natural action of the wreath product $\Gn = \G^n
\rtimes S_n$. The direct sum of equivariant $K$-groups
$\bigoplus_{n \ge 0}K_{\Gn}(X^n)  \otimes \C$ were shown earlier
by the author to carry several interesting algebraic structures.
In this paper we study the $K$-groups $K_{\tG_n}(X^n)$ of
$\Gn$-equivariant Clifford supermodules on $X^n$. We show that
$\tFG =\bigoplus_{n\ge 0}K_{\tG_n}(X^n) \otimes \C$ is a Hopf
algebra and it is isomorphic to the Fock space of a twisted
Heisenberg algebra. Twisted vertex operators make a natural
appearance. The algebraic structures on $\tFG$ , when $\G$ is
trivial and $X$ is a point, specialize to those on a ring of
symmetric functions with the Schur $Q$-functions as a linear
basis. As a by-product, we present a novel construction of
$K$-theory operations using the spin representations of the
hyperoctahedral groups.
 }
\end{abstract}
\tableofcontents
\section{Introduction}
Motivated in part by Vafa-Witten \cite{VW1} and generalizing the
work of Segal \cite{Seg2} (also cf. Grojnowski \cite{Gr}), we
studied in \cite{W1} a direct sum, denoted by $\FG$, of the
equivariant $K$-groups $\oplus_{n\geq 0} K_{\Gn}(X^n) \otimes \C$
associated to a topological $\G$-space $X$. Here $\G$ is a finite
group and the wreath product $\Gn$ (i.e. the semi-direct product
$\G^n \rtimes S_n$) acts naturally on the $n$th Cartisian product
$X^n$. We proved that the space $\FG$ carries several remarkable
algebraic structures such as Hopf algebra and Fock representation
of a Heisenberg (super)algebra etc, and that vertex operators
makes a natural appearance as a part of $\la$-ring structure. We
in addition pointed out in \cite{W1} a new approach to the
realization of the Frenkel-Kac-Segal homogeneous vertex
representations of affine Lie algebras by using representation
rings of the wreath product $\Gn$ associated to a finite subgroup
$\G$ of $SL_2 (\C)$. This has been subsequently completed in
\cite{FJW1} jointly with I.~Frenkel and Jing, and further extended
in \cite{FJW2} to realize vertex representations of twisted affine
and toroidal Lie algebras by using the spin representation rings
of a double cover $\widetilde{\G}_n$ of the wreath product $\Gn$.

In this paper we will introduce a variant of equivariant
$K$-theory. Given a topological space $X$ acted upon by a finite
supergroup $G$ in an appropriate sense, we introduce a category of
complex $G$-equivariant {\em spin} vector {\em super}bundles over
$X$, and consider the corresponding Grothendieck group $K^-_G(X)$.
The superscript $-$ here and below is used in this paper to stand
for {\em spin}, i.e. a certain central element $z$ in the
supergroup $G$ acts as $-1$. This $K$-group can be thought as an
invariant of orbifolds with (certain distinguished) discrete
torsion as introduced by Vafa \cite{Va}. Discrete torsion has been
a topic of interest from various viewpoints since then, cf.
\cite{VW2, Di, AR, Ru, Sh} and the references therein. We present
here a variant of the decomposition theorem of Adem and Ruan
\cite{AR} for what they call twisted equivariant $K$-theory, which
generalizes the decomposition theorem \cite{BC, Ku} (also cf.
\cite{AS, HKR}) in equivariant $K$-theory.

Our formulation of such $K$-groups is motivated by providing a
general framework for the main subject of study in this paper,
namely a spin/twisted version of the space $\FG$ studied in
\cite{W1}. When the topological space under consideration is a
point, our $K$-theory specializes to the theory of supermodules of
finite supergroups, cf. J\'ozefiak \cite{Jo1}. (In this paper {\em
super} always means {\em $\Z_2$-graded}). Such a theory of
supermodules has provided, in our opinion, a most natural
framework (cf. \cite{Jo}) for the exposition of the spin
representations of a double cover $\widetilde{S}_n$ of the
symmetric group initiated by Schur \cite{Sc}. The spin
representation theory of a double cover $\widetilde{H}_n$ of the
hyperoctahedral group, being almost parallel to and to some extent
simpler than the spin representation theory for $\widetilde{S}_n$,
can also be treated successfully in terms of supermodules (cf.
\cite{Jo2}).

Given a space $X$ acted by a (non-graded) finite group $\G$, we
obtain our main example of such $K$-group $K^-_{\tG_n}(X^n)$ by
considering the action on the $n$-th Cartesian product $X^n$ by
the wreath product $\Gn$ which is further extended trivially to
the action of a larger finite supergroup $\tG_n$. Here $\tG_n$ is
a double cover of the wreath product $(\G \times \Z_2)^n \rtimes
S_n$. (We recommend that the reader sets $\G$ to be the
one-element group in their first reading so that the whole picture
becomes simpler and more transparent. In this case $\tG_n$ reduces
to a double cover $\widetilde{H}_n$ of the Hyperoctahedral group.)
The category of $\tG_n$-equivariant spin vector superbundles over
$X^n$ admits an equivalent reformulation which has a perhaps
better geometric meaning. Namely this is the category of
$\Gn$-equivariant vector bundles $E$ over $X^n$ such that $E$
carries a supermodule structure with respect to the complex
Clifford algebra of rank $n$ which is compatible with the action
of $\Gn$.

A fundamental example of $\tG_n$-vector superbundles over $X^n$,
which plays an important role in this paper, is given as follows
for $X$ compact. Given a $\G$-vector bundle $V$ over $X$, we
consider the vector superbundle $V|V =V\oplus V $ over $X$ with
the natural $\Z_2$-grading. We can endow the $n$th outer tensor
$(V|V)^{\boxtimes n}$ a natural $\tG_n$-equivariant vector
superbundle structure over $X^n$.

We will show that the direct sum $$
 \tFG := \bigoplus_{n=0}^{\infty} K^-_{\tG_n}(X^n)\bigotimes \C
 $$
carries naturally a Hopf algebra structure and it is isomorphic to
the Fock space of a twisted Heisenberg superalgebra associated to
$K^-_{\tG_1}(X) \cong K_\G(X)$. All the algebraic structures are
constructed in terms of natural $K$-theory maps. In particular the
dimension of $K^-_{\tG_n}(X^n)$ is determined explicitly for all
$n$. We remark that such a twisted Heisenberg algebra has played
an important role in the theory of affine Kac-Moody algebras, cf.
\cite{FLM}. Roughly speaking, the structure of the space $\FG$
studied in \cite{Seg2, W1} is modeled on the ring $\Lambda_{\C}$
of symmetric functions with a basis given by Schur functions (or
equivalently on the direct sum of representation rings of
symmetric groups $S_n$ for all $n$). The structure of the space
$\tFG$ under consideration of this paper is shown to be modeled
instead on the ring $\Omega_{\C}$ of symmetric functions with a
linear basis given by the so-called Schur $Q$-functions (or
equivalently on the direct sum of the spin representation ring of
$\Hy_n$ for all $n$). It is well known that the graded dimension
of the ring $\Omega_{\C}$ is given by the generating function
 $$
  \frac1{\prod_{i=0}^\infty (1
 -q^{2i-1})} =
 q^{-\frac{1}{24}}\frac{\eta(q^2)}{\eta(q)},
 $$
where $\eta(q)$ is the Dedekind $\eta$ function.

Just as the $\la$ ring is modeled on the ring $\Lambda_{\C}$ of
symmetric functions (cf. e.g. \cite{Kn}), one can introduce a
$Q$-$\la$ ring\footnote{Here $Q$ stands for {\em Queer} or Schur
$Q$-functions. We believe that there exists a rich
$Q$-mathematical world which is relevant to various twisted, spin,
super structures, etc.} structure modeled on the ring
$\Omega_{\C}$ with Adams operations of odd degrees only. We show
that as a part of the $Q$-$\la$ ring structure on $\tFG$ twisted
vertex operators naturally appear in $\tFG$ via the $n$-th outer
tensor $(V|V)^{\boxtimes n}$ associated to $V\in K_{\G}(X)$ in
terms of the Adams operations.

It is of independent interest to see that when we restrict
$K^-_{\Hy_n}(X^n)$ from $X^n$ to its diagonal we are able to
obtain various $K$-theory operations on $K(X)$, including
supersymmetric power operations and Adams operations of odd
degrees, by means of the spin supermodules of $\Hy_n$. This is a
super analog of Atiyah's construction \cite{At} of $K$-theory
operations on $K(X)$ by means of the representations of the
symmetric groups.

Motivated by G\"ottsche's formula, Vafa and Witten \cite{VW1}
conjectured that the direct sum ${\cal H}(S)$ of homology groups
for Hilbert schemes $S^{[n]}$ of $n$ points on a
(quasi-)projective surface $S$ should carry the structure of a
Fock space of a Heisenberg algebra, which was realized
subsequently in a geometric way by Nakajima and Grojnowski (cf.
\cite{Na, Gr}). Parallel algebraic structures such as Hopf
algebra, vertex operators, and Heisenberg algebra as part of
vertex algebra structures \cite{Bo, FLM}, have naturally showed up
in ${\cal H}(S)$ as well as in $\FG$. When $S$ is a suitable
resolution of singularities of an orbifold $X/\G$, there appears
close connections between ${\cal H}(S)$ and $\FG$, cf. \cite{W2,
W3} and the references therein. It will be very important to see
if one can find a `Hilbert scheme' version of the orbifold picture
drawn in this paper. It is also interesting to see if our current
work can find some applications in string theory, cf. \cite{VW1,
VW2, Di, Sh}. In fact the special case of our construction for
$\G$ trivial is closely related to an earlier paper of Dijkgraaf
\cite{Di}. See the Appendix.

When $X$ is a point, the $K$-group $\underline{K}^-_{\tG_n}(X^n)$
becomes the Grothendieck group of spin supermodules of $\tG_n$. In
a companion paper \cite{JW} joint with Jing, the conjugacy classes
and Grothendieck groups of $\tG_n$ have been studied in detail. In
particular, when $\G$ is a subgroup of $SL_2(\C)$, they are used
to realize vertex representations of twisted affine algebras (cf.
\cite{FLM}) and toroidal Lie algebras. This provides a new group
theoretic construction, a somewhat improved version in our opinion
than the one in \cite{FJW2} using the groups $\widetilde{\G}_n$,
of twisted vertex representations. In the Appendix, we sketch
another formulation of our main results in this paper using the
groups $\widetilde{S}_n$ and $\widetilde{\Gn}$ instead of
$\widetilde{H}_n$ and $\tG_n$.

The plan of this paper is as follows. In
Sect.~\ref{sec_supergroup} we present the spin representation
rings of the finite supergroup $\tG_n$. In Sect.~\ref{sec_equiv}
we introduce the category and $K$-groups of spin vector
superbundles which are equivariant with respect to finite
supergroups and present a decomposition theorem for such
equivariant $K$-groups. In Sect.~\ref{sec_oper} we study
$K$-theory operations based on the spin representations of
$\Hy_n$. The results in this section are not to be further used in
this paper. In Sect.~\ref{sec_structure}, which is the heart of
the paper, we present the structures of a Hopf algebra and of a
$Q$-$\la$ ring on $\tFG$, and relate the latter to the twisted
vertex operators. We further construct a Heisenberg algebra which
acts on $\tFG$ irreducibly by means of natural $K$-theory maps. In
the Appendix, we outline a somewhat different construction in
terms of the group $\widetilde{\Gn}$.
\section{The group $\bgn$ and its spin supermodules}
\label{sec_supergroup}

In this section we recall briefly some essential points in the
theory of supermodules of a finite supergroup (cf. \cite{Jo1}). We
define the finite supergroup $\tG_n$ associated to any finite
group $\G$, and study its conjugacy classes and spin supermodules.
More detail of these can be found in \cite{JW}.
\subsection{Definition of the supergroup $\tG_n$}
Let $\GG$ be a finite group and let $d: \GG \rightarrow \Z_2$ be a
group epimorphism. We denote by $\GG_0$ the kernel of $d$ which is
a subgroup of $\GG$ of index $2$. We regard $d(\cdot)$ as a parity
function on $\GG$ by letting the degree of elements in $\GG_0$ be
$0$ and the degree of elements in the complementary $\GG_1 = \GG
\backslash \GG_0$ be $1$. Elements in $\GG_0$ (resp. $\GG_1$) will
be called even (resp. odd). We will often refer to the pair $(\GG,
d)$, or simply $\GG$ when there is no ambiguity, as a {\em finite
supergroup}. The class of finite supergroups under consideration
in this paper has an additional property: it contains a
distinguished even central element $z$ of order $2$. We denote by
$\theta$ the quotient group homomorphism $\GG \rightarrow G \equiv
\GG/ \langle 1, z \rangle.$

In this paper the modules over a finite supergroup or a
superalgebra (such as the group superalgebra of a finite
supergroup) will always be $\Z_2$-graded (i.e. supermodules)
unless otherwise specified. A general theory of supermodules over
finite supergroups was developed by J\'ozefiak \cite{Jo1}. This
was motivated to provide a modern account \cite{Jo} of Schur's
seminal work on spin representations of symmetric groups
\cite{Sc}.

Given two supermodules $M=M_0 +M_1$ and $N=N_0 +N_1$ over a
superalgebra $A =A_0 +A_1$, the linear map $f: M\rightarrow N$
between two $A$-supermodules is a homomorphism of degree $i$ if
$f(M_j)\subset M_{i+j}$ and for any homogeneous element $a\in A$
and any homogeneous vector $m\in M$ we have
\[ f(am)=(-1)^{d(f)d(a)}a f(m).
\]
The degree $0$ (resp. $1$) part of a superspace is referred to as
the {\em even} (resp. {\em odd}) part. We denote
\[
\Hom_A(M, N)=\Hom_A(M, N)_0\oplus \Hom_A(M, N)_1,
\]
where $\Hom_A(M, N)_i$ consists of $A$-homomorphisms of degree $i$
from $M$ to $N$. The notions of submodules, tensor product, and
irreducibility etc for supermodules are defined similarly.

Given a finite supergroup $\GG$, a $\GG$-supermodule $V$ is called
{\em spin} if the central element $z$ acts as $-1$. Alternatively,
one can associate a $2$-cocycle $\alpha: G \times G
\rightarrow\Z_2 =\{\pm 1\}$, such that $V$ becomes a projective
supermodule of the group $G =  \GG/ \langle 1,z \rangle$
associated with $\alpha$, namely, the actions of any two elements
$g, h \in G$ on $V$, denoted by $\rho(g), \rho(h)$, satisfy the
relation
\begin{eqnarray} \label{eq_cycle}
  \rho(g) \rho(h) = \alpha(g,h)\rho (gh).
\end{eqnarray}
Among all $\GG$-supermodules, we will only consider the spin
supermodules in this paper. It is clear that the restriction
(resp. the induction) of a spin supermodule to a $\Z_2$-graded
subgroup (resp. a larger supergroup) with the same distinguished
even element $z$ remains to be a spin supermodule.

There are two types of complex simple superalgebras according to
C.T.C.~Wall: $M(r|s)$ and $Q(n)$. Here $M(r|s)$ is the
superalgebra consisting of the linear transformations of the
superspace $\C^{r|s} =\C^r + \C^s$. The superalgebra $Q(n)$ is the
graded subalgebra of $M(n|n)$ consisting of matrices of the form
\begin{eqnarray*}
\left[
\begin{array}{cc} A & B\\ B & A
\end{array}
\right] .
\end{eqnarray*}
It is known \cite{Jo1} that the group (super)algebra of a finite
supergroup is semisimple, i.e. decomposes into a direct sum of
simple superalgebras. According to the classification of simple
superalgebras above, the irreducible supermodules of a finite
supergroup are divided into two types, type $M$ and type $Q$. We
note that the endomorphism algebra of an irreducible supermodule
$V$ is isomorphic to $\C$ if $V$ is of type $M$ and isomorphic to
the complex Clifford algebra $C_1$ in one variable if $V$ is of
type $Q$.

Let $\NP_n$ the group generated by $1,z,a_1,\ldots,a_n$ subject to
the relations
\begin{eqnarray*}
z^2=1,\ a_i^2=z\ {\rm and}\  a_ia_j=z a_ja_i,\ {\rm for}\ i\not=j.
\end{eqnarray*}
The symmetric group $S_n$ acts on $\NP_n$ via permutations of the
elements $a_1,\ldots,a_n$, i.e. $\sigma (a_i) = a_{\sigma (i)}$
for $\sigma \in S_n$. We thus form the semi-direct product $\Hy_n
=\NP_n \rtimes S_n$, which naturally endows a $\Z_2$-grading given
by the parity $d(a_i)=1$, $d(z) =0$, and $d(\sigma)=0$ for
$\sigma\in S_n$. We may therefore regard $\Hy_n$ as a finite
supergroup with a distinguished even central element. Note that
the group superalgebra $\C[\Hy_n] / \langle z =-1 \rangle$ is
exactly the complex Clifford algebra $C_n$ in $n$ variables.

Let $\G$ be a finite group with $r +1$ conjugacy classes. We
denote by $\G^*=\{\g_i\}_{i=0}^{r}$ the set of complex irreducible
characters where $\g_0$ denotes the trivial character, and by
$\G_*$ the set of conjugacy classes. Let $R(\G)=\oplus_{i=0}^{r}
\C\g_i$ be the space of class functions on $\G$, and set
$R_{\Z}(\G)=\oplus_{i=0}^{r} \C\g_i$. For $c \in \G_*$ let
$\zeta_c$ be the order of the centralizer of an element in the
conjugacy class $c$, so the order of the class is then $|\G
|/\zeta_c$.

Given a positive integer $n$, let $\G^n = \G \times \cdots \times
\G$ be the $n$-th direct product of $\G$, and let $\G^0$ be the
trivial group. The symmetric group $S_n$ naturally acts on $\G^n
\times \NP_n $ by simultaneous permutations of elements in $\G^n$
and $ \NP_n $.

The finite supergroup $\tG_n$ is then defined to be the
semi-direct product of the symmetric group $S_n$ and $\G^n \times
\NP_n $, with the multiplication given by
$$
 (g, \sigma)\cdot (h, \tau)=(g \, \sigma(h), \sigma \tau),
 \quad g,h \in \G^n \times \NP_n, \; \sigma, \tau \in S_n.
 $$
The order of $\tG_n$ is clearly $2^{n+1} n!|\G|^n$. The
$\Z_2$-grading on $\tG_n$ is induced from that on $\NP_n$ and by
letting the elements in $\G^n$ be even (i.e. of degree $0$).
Denoting by $\bG_n = (\G \times \Z_2)^n \rtimes S_n$, we have the
following exact sequence of groups:
 $$
 1 \rightarrow \Z_2 = \{1, z\} \rightarrow \tG_n
 \stackrel{\theta_n}{\rightarrow} \bG_n
 \rightarrow 1.
$$ It is clear that when $\G$ is trivial $\tG_n$ reduces to a
double cover $\Hy_n$ of the hyperoctahedral group $H_n = \Z_2^n
\rtimes S_n$.

The finite supergroup $\tG_n$ contains $\NP_n, \Hy_n$ and the
wreath product $\Gn =\G^n \rtimes S_n$ as distinguished subgroups.
Letting $\G^n$ act trivially on $\NP_n$ we extend the action of
the symmetric group $S_n$ to $\Gn$ on $\NP_n$. In this way we may
also view $\tG_n$ as a semi-direct product between $\Gn$ and
$\NP_n$.
\subsection{Conjugacy classes of $\tG_n$}
Let $\GG$ be a finite supergroup and  put $G = \GG/ \langle
1,z\rangle$ and $\theta: \GG \rightarrow G$ as before. For any
conjugacy class $C$ of $G$, $\theta^{-1} (C)$ is either a
conjugacy class of $\GG$ or it splits into two conjugacy classes
of $\GG$, cf. \cite{Jo1}. If $\theta^{-1} (C)$ splits, the
conjugacy class $C$ will be referred to as {\em split}, and an
element $g$ in $C$ is also called {\em split}, which is equivalent
to say that the two preimages of $g$ under $\theta$ are not
conjugate to each other. Otherwise $g$ is said to be {\it
non-split}. In view of (\ref{eq_cycle}), we have the following
easy equivalent formulation.

\begin{lemma} \label{lem_chara}
An element $g$ in $G = \GG/ \langle 1,z\rangle$ is split if and
only if
 $$\varepsilon_g(\cdot):= \alpha (g, \cdot)\alpha (\cdot,g)^{-1}
 $$
defines a {\em trivial} character of the centralizer group
$Z_G(g)$.
\end{lemma}

For the study of spin supermodules of $\tG_n$, it is crucial to
have a detailed description of split conjugacy classes of $\GG$.
Indeed the characters of spin supermodules vanish on non-split
conjugacy classes as well as on odd split classes, cf. \cite{Jo1}.
Below we will concentrate on the group $\tG_n$.

Let $\la=(\la_1, \la_2, \cdots, \la_l)$ be a partition of integer
$|\la|=\la_1+\cdots+\la_l$, where $\la_1\geq \dots \geq \la_l \geq
1$. The integer $l$ is called the {\em length} of the partition
$\la $ and is denoted by $l (\la )$. We will also make use of
another notation for partitions: $$ \la=(1^{m_1}2^{m_2}\cdots) ,
$$ where $m_i$ is the number of parts in $\la$ equal to $i$. A
partition $\la$ is {\it strict} if its parts are distinct
integers, namely all the multiplicities $m_i$ are $1$ or $0$.
Given a partition $\lambda = (1^{m_1} 2^{m_2} \ldots )$ of $n$, we
define
\[
  z_{\la } = \prod_{i\geq 1}i^{m_i}m_i!.
\]
We note that $z_{\la }$ is the order of the centralizer of an
element of cycle-type $\la $ in $S_n$.

For a finite set $X$ and $\rho=(\rho(x))_{x\in X}$ a family of
partitions indexed by $X$, we write $$\|\rho\|=\sum_{x\in
X}|\rho(x)|.$$ It is convenient to regard $\rho=(\rho(x))_{x\in
X}$ as a partition-valued function on $X$. We denote by ${\mathcal
P}(X)$ the set of all partitions indexed by $X$ and by ${\mathcal
P}_n(X)$ the set of all partitions in ${\mathcal P}(X)$ such that
$\|\rho\|=n$. The total number of parts, denoted by
$l(\rho)=\sum_xl(\rho(x))$, is called the {\em length} of $\rho$.
Let ${\cal OP}(X)$ be the set of partition-valued functions
$(\rho(x))_{x\in X}$ in ${\mathcal P}(X)$ such that all parts of
the partitions $\rho(x)$ are odd integers, and let ${\mathcal
SP}(X)$ be the set of partition-valued functions $\rho\in
{\mathcal P}(X)$ such that each partition $\rho(x)$ is strict. It
is clear that $|{\cal OP}_n(X)|=|{\cal SP}_n(X)|. $ When $X$
consists of a single element, we will omit $X$ and simply write
$\mathcal P$ for ${\mathcal P}(X)$, thus $\mathcal OP$ or
$\mathcal SP$ will be used similarly.

We denote by
\begin{eqnarray*}
{\cal P}_n^+(X)&=\{\la\in {\cal P}_n(X)|\quad l(\rho)\equiv 0
\bmod{2} \}, \\
 {\cal P}_n^-(X) &=\{\la\in {\cal P}_n(X)|\quad l(\rho)\equiv 1
\bmod{2} \},
\end{eqnarray*}
and define ${\cal SP}_n^{\pm}(X)={\mathcal P}_n^{\pm}(X)\cap
{\mathcal SP}_n(X)$ for $i=0, 1$.

The conjugacy classes of a wreath product is well understood, cf.
\cite{M}. In particular this gives us the following description of
conjugacy classes of the wreath product $\bG_n =(\G \times \Z_2)^n
\rtimes S_n$.

Let $x=(g, \sigma)$ be an element in a conjugacy class of $\tG_n$,
where $g=(g_1, \cdots, g_n)$ and $g_i = (\alpha_i, \varepsilon_i)
\in\G\times \Z_2$. We take the convention here that $\Z_2 =\{ \pm
1 \}$. For each cycle $y=(i_1 i_2 \cdots i_k)$ in the permutation
$\sigma$ consider the element $\alpha_y=\alpha_{i_k} \alpha_{i_{k
-1}} \cdots \alpha_{i_1} \in \G$ and $\varepsilon_y =
\varepsilon_{i_k} \varepsilon_{i_{k -1}} \cdots \varepsilon_{i_1}$
(which is $\pm 1$) corresponding to the cycle $y$. For each
$c\in\G_*$, $\varepsilon = \pm$ and $r\geq 1$, let
$m^{\varepsilon}_r(c)$ be the number of $r$-cycles in the
permutation $\sigma$ such that the cycle product $\alpha_y$ lie in
the conjugacy class $c$ and $\varepsilon_y$ equals $\varepsilon
\cdot 1$. Then $c\rightarrow
\rho_{\varepsilon}(c)=(1^{m^{\varepsilon}_1(c)}
2^{m^{\varepsilon}_2(c)}\ldots)$ defines a partition-valued
function on $\G_*$ for each $\varepsilon$. The partition-valued
function $\rho_{\varepsilon} = (\rho_{\varepsilon}(c))_{c\in\G_*}
\in {\cal P}(\G_*)$ is called the {\em $\varepsilon$-type} of $x$.
Denote by ${\cal P}(\G_*)^2_n$ the set of pairs of
partition-valued functions $\rho=(\rho^+, \rho^-)$ such that
$||\rho^+|| + ||\rho^-|| =n$. The pair $\rho=(\rho^+, \rho^-)$ is
called the {\em type} of $x$. One can show that the type only
depends on the conjugacy class of $x$ in $\bG_n$ and the conjugacy
classes in $\bG_n$ are parameterized by the types $\rho \in {\cal
P}(\G_*)^2_n$. We will also say that the conjugacy class
containing $x$ has conjugacy type $\rho$ and is denoted by
$C_{\rho}$ if $x$ is of type $\rho$.

Denote by $D_{\rho} = \theta_n^{-1} (C_{\rho})$. The following is
established in \cite{JW}, Theorem 2.1.

\begin{theorem} \label{th_splitclass}
For $\rho = (\rho^+, \rho^-) \in {\cal P}(\G_*)^2_n$, $D_{\rho}$
splits into two conjugacy classes in $\tG_n$ if and only if:

(1) for $D_{\rho}$ is even we have $\rho^+ \in {\cal OP}_n(\G_*)$
and $\rho^- = \emptyset$;

(2) for $D_{\rho}$ is odd we have $\rho^+ = \emptyset$ and $\rho^-
\in {\cal SP}_n^-(\G_*)$.
\end{theorem}

Thus, the set $(\bG_n)_*^{e.s}$ of even split conjugacy classes in
$\bG_n$ is in one-to-one correspondence with the set ${\cal
OP}_n(\G_*)$.
\subsection{Spin supermodules over $\tG_n$}

The number of irreducible spin supermodules of $\tG_n$ equals the
number of even split conjugacy classes in $\tG_n$, by a general
theorem in the supermodule theory of finite supergroups
\cite{Jo1}. The next proposition follows from the equality $|{\cal
SP}_n(\G_*)| =|{\cal OP}_n(\G_*)|$.

\begin{proposition}
The number of irreducible spin supermodules of $\tG_n$ is equal to
the number $|{\cal SP}_n(\G_*)|$ of strict partition-valued
functions on $\G_*$.
\end{proposition}

We denote by $R^-(\tG_n)$ (resp. $R^-_{\Z}(\tG_n)$) the $\C$-span
(resp. $\Z$-span) of the characters of irreducible spin
supermodules of $\tG_n$. Let
$$
 {R}^-_{\G} = \bigoplus_{n=0}^{\infty} R^-(\tG_n).
$$
 When $\G$ is trivial, we will simply drop the subscript $\G$
and write
 $$R^- =\bigoplus_{n=0}^{\infty} R^-(\Hy_n).$$

For example, when $\G$ is trivial and thus $\tG_n$ reduces to
$\Hy_n$, the irreducible spin supermodules of $\Hy_n$ are
parameterized by strict partitions of $n$ (cf. \cite{Ser} and
\cite{Jo2}). For strict partitions $\la$ and $\mu$ of $n$, let
$T^\la$ and $T^\mu$ denote the corresponding irreducible spin
supermodules over $\Hy_n$. We have

\begin{equation}\label{octahedral2}
{\rm dim }\,{\rm
Hom}_{\Hy_n}(T^\la,T^\mu)=\delta_{\la\mu}2^{\delta(l(\la))},
\end{equation}
where the number $\delta(l(\la))$ is $0$ for $l(\la)$ even and $1$
otherwise. That is, the supermodule $T^{\la}$ is of type $M$
(resp. type $Q$) if and only if $l(\la)$ is even (resp. odd).

A most distinguished example of irreducible $\Hy_n$-supermodule is
the so-called {\em basic spin supermodule} $L_n$ constructed as
follows (cf. \cite{Jo2, JW}). As a superspace $L_n$ is isomorphic
to the group superalgebra $\C[\NP_n]/\langle z=-1\rangle$ (i.e.
the Clifford algebra in $n$ variables). Denote by $y_i \in L_n$
the image of $a_i \in \NP_n$. Then $y_I =\prod_{i \in I} y_i, I
\subset \{1, \ldots, n\}$, form a linear basis of $L_n$. The
action of $\Hy_n$ on $L_n$ is given by
\[
a_j y_I =y_j y_I \;(j =1, \ldots, n), \quad \sigma y_I = y_{s(I)},
\; s \in S_n.
\]
Indeed $L_n$ is exactly the $\Hy_n$-supermodule $T^{(n)}$
associated to the one-part partition $(n)$. If we denote by
$\xi^n$ the character of $L_n$, then the character value of
$\xi^n$ is $2^{l(\rho)}$ on a conjugacy class of type $\rho
=(\rho^+, \emptyset)$, where $\rho^+ \in {\cal OP}_n(\G_*)$.

For each partition-valued function $\rho=(\rho(c))_{c\in\G_*}$ we
define $$ Z_{\rho}=2^{l(\rho)}\prod_{c\in\G_*}z_{\rho(c)}
\zeta_c^{l(\rho(c))}, $$ which is the order of the centralizer of
an element in $\bG_n$ of conjugacy type $\rho=( \rho(c)
)_{c\in\G_*}$ (see \cite{JW}).

For a fixed $c \in \G_*$, we denote by $c_n$ ($n \in 2\Z_+ +1$)
the even split conjugacy class in $\bG_n$ of the type $(\rho^+,
\emptyset)$, where the partition-valued function $\rho^+$ takes
value the one-part partition $(n)$ at $c \in \G_*$ and zero
elsewhere. We denote by $\sigma_n (c) \in R^-(\tG_n)$ the class
function of $\tG_n$ which takes value $n \zeta_c$ at the conjugacy
class $c_n$ $(c \in \G_*)$ and zero otherwise. For $\rho = \{m_r
(c) \}_{c,r} \in  {\cal OP}_n(\G_*)$, we define $$\sigma^{\rho} =
\prod_{c \in G_*, r \geq 1} \sigma_r(c)^{m_r(c)}.$$ and regard it
as the class function on $\tG_n$ which takes value $Z_{\rho}$ at
the conjugacy class $D^+_{\rho}$ and $0$ elsewhere. Then it
follows that (cf. \cite{Ser, Jo2})

\begin{eqnarray} \label{eq_basicspin}
 \xi^n = \sum_{\rho\in{\cal OP}_n(\G_*)}
      2^{l(\rho)} Z_{\rho}^{ -1}\sigma^{\rho}.
\end{eqnarray}

Finally we define the analog for $\tG_n$ of Young subgroups of the
symmetric groups. Let $\tG_n\tilde{\times}\tG_m$ be the direct
product of $\tG_n$ and $\tG_m$ with a twisted multiplication
\begin{eqnarray*}
(t, t')\cdot (s, s')=(ts z^{d(t')d(s)}, t's'),
\end{eqnarray*}
where $s, t\in\tG_n, s', t'\in\tG_m$. We define the {\it spin
product} of $\tG_n$ and $\tG_m$ by letting
\begin{eqnarray} \label{eq_product}
\tG_n\hat{\times}\tG_m={\tG_n\tilde{\times}\tG_m}/ \{(1, 1), (z,
z)\}.
\end{eqnarray}
which carries a canonical $\Z_2$-grading and can be regarded as a
subgroup of the supergroup $\tG_{n+m}$ in a natural way.

For two spin supermodules $U$ and $V$ of $\tG_n$ and $\tG_m$ we
define a $\tG_n\hat{\times}\tG_m$ spin supermodule on the tensor
product $U\otimes V$ by letting
\begin{eqnarray*}
(t, s)\cdot (u\otimes v)=(-1)^{d(s)d(u)}(tu\otimes sv).
\end{eqnarray*}
This induces an isomorphism $\phi_{n, m}$ from $R^-(\tG_n )
\bigotimes R^-(\tG_m)$ to $R^-(\tG_n \hat{\times} \tG_m)$.

The space $\tRG$ carries a (commutative associative)
multiplication which is defined by the composition (for all $n,
m$)
\begin{eqnarray*}
   R^-(\tG_n ) \bigotimes R^-(\tG_m)
 \stackrel{\phi_{n, m} }{\longrightarrow} R^-(\tG_n \hat{\times} \tG_m)
 \stackrel{Ind}{\longrightarrow} R^-( \tG_{n + m}).
\end{eqnarray*}
\section{A decomposition theorem in equivariant $K$-theory} \label{sec_equiv}
In this section we introduce a variant of equivariant $K$-theory
for a finite supergroup. We recall a decomposition theorem in the
equivariant $K$-theory from \cite{BC, Ku, AS, HKR}, and present a
generalization of it in our new setup.
\subsection{The standard version}
Given a (non-graded) finite group $\G$ and a compact Hausdorff
$\G$-space $X$, we recall \cite{Seg1} that $K^0_\G(X)$ is the
Grothendieck group of $G$-vector bundles over $X$. One can define
$K_\G^1 (X)$ in terms of the $K^0$ functor and a certain
suspension operation, and one puts
 $$
 K_\G(X) = K^0_\G(X)\bigoplus K^1_\G(X).
 $$
 In this paper we will be only concerned
about $K_\G(X) \otimes \C$, and subsequently we will denote
 $$
 \underline{K}_\G(X) = K_\G(X) \bigotimes \C.
 $$
We denote by $\dim {K}_\G^i (X) (i = 0, 1)$ the dimension of
${K}_\G^i (X) \otimes \C$.

If $X$ is locally compact, Hausdorff and paracompact $\G$-space,
take the one-point compactification $X^+$ with the extra point
$\infty$ fixed by $\G$. Then we define $K_\G^0 (X)$ to be the
kernel of the map $$ K^0_\G(X^+) \longrightarrow K^0_\G(
\{\infty\} ) $$ induced by the inclusion map $\{\infty\}
\hookrightarrow X^+$. This definition is equivalent to the earlier
one when $X$ is compact. We also define $K^1_\G(X) = K^1_\G(X^+).$

Note that $K_{\G} (pt)$ is isomorphic to the Grothendieck ring
$R_{\Z}(\G)$ and $\underline{K}_{\G} (pt)$ is isomorphic to the
ring $R(\G)$ of class functions on $\G$.

Let $X^g$ denote the fixed-point set by $g \in \G$, which is
preserved under the action of the centralizer $Z_\G(g)$. The
following decomposition theorem (cf. \cite{BC, Ku, AS, HKR}) gives
a description of each direct summand over conjugacy classes of
$\G$. We remark that the subspace of invariants $\underline{K}
(X^{g})^{Z_{\G}(g)}$ is isomorphic to $\underline{K}(X^{g}
/Z_{\G}(g))$, and it is isomorphic for different choice of $g$ in
the same conjugacy class $[g] \in \G_*$.

\begin{theorem} \label{th_tech}
  There is a natural $\Z_2$-graded isomorphism
  \begin{eqnarray*}
   \phi : \underline{K}_{\G} (X)  \longrightarrow
    \bigoplus_{[g]\in \G_*} \underline{K} \left( X^{g}
                                  \right)^{Z_{\G} (g)}.
  \end{eqnarray*}
\end{theorem}
\subsection{A super/twisted variant}
\label{subsect_ktwist}

Now let $\GG$ be a finite supergroup which contains a
distinguished central element $z$ of order $2$, and let $\theta$
be the quotient homomorphism $\GG \rightarrow  {G} = \GG/\langle
1,z \rangle.$

Let $X$ be a compact (non-graded) ${G}$-space. We may regard $X$
as a $\GG$-space by letting $g \in \GG$ act by $\theta (g) \in G$.
We denote by ${\cal C}^-_{\GG}(X)$ the category whose objects
consist of $\GG$-equivariant complex vector superbundle (i.e.
$\Z_2$-graded bundles) $E = E_0 +E_1$ (often denoted by $E_0|E_1$)
over $X$ on which $\GG$ acts in a $\Z_2$-graded manner and $z$
acts as $-1$. We will call such an object a  {\em spin}
$\GG$-vector superbundle. Given two objects $E,F$ in the category
${\cal C}^-_{\GG}(X)$, the space of homomorphisms of
$\GG$-equivariant vector superbundles between $E$ and $F$ admits a
natural $\Z_2$-grading:

 $$\Hom^{\GG} (E, F) = \Hom^{\GG}_0
(E, F) \bigoplus \Hom^{\GG}_1 (E, F).
 $$
By our definition of finite supergroups, $\GG$ always contains odd
elements. It follows that $\rk E_0 = \rk E_1$ thanks to the
existence of odd automorphisms from the odd elements of $\GG$.

We denote by ${K}^{-,0}_{\GG}(X)$ the Grothendieck group of the
abelian monoid of isomorphism classes\footnote{The isomorphisms
are $\Z_2$-graded and isomorphisms of degree $1$ are allowed.} of
the vector superbundles in ${\cal C}^-_{\GG}(X)$. As in the
ordinary case, we can extend the definition of ${\cal
C}^-_{\GG}(X)$ to locally compact spaces, and define
${K}^{-,1}_{\GG}(X)$ to be ${K}^{-,0}_{\GG}(X \times {\Bbb R})$
where $\Bbb R$ is the real line. We denote ${K}^-_{\GG}(X) =
{K}^{-,0}_{\GG}(X) + {K}^{-,1}_{\GG}(X)$.

In this paper we will be only concerned about the free part
${K}^-_{\GG}(X) \otimes \C$, which will be denoted by
$\underline{K}^-_{\GG}(X)$ subsequently. The following theorem,
generalizing Theorem~\ref{th_tech}, is a variation of Adem-Ruan's
construction \cite{AR} for the so-called twisted equivariant
$K$-theory. Recall that the character $\varepsilon_g$ of the
centralizer group $Z_G(g)$ was defined in Lemma~\ref{lem_chara}
and $Z_G(g)$ acts on $\underline{K}(X^{g})$.

\begin{theorem} \label{th_supertech}
Let $\GG$ be a finite supergroup which contains a distinguished
central element $z$ of order $2$, and let $\theta$ be the quotient
homomorphism $\GG \rightarrow G =\GG/\langle 1, z \rangle$. Given
a locally compact Hausdorff $G$-space $X$ and regarding it as a
$\GG$-space, we have a natural $\Z_2$-graded isomorphism
  \begin{eqnarray} \label{eq_decompo}
   \phi : \underline{K}^-_{\GG}(X)  \stackrel{\cong}{\longrightarrow}
    \bigoplus_{[g]}
     \left( \underline{K}(X^{g}) \bigotimes  \varepsilon_g
     \right)^{Z_{G} (g)},
  \end{eqnarray}
where the summation runs over the even conjugacy classes in $G$.
\end{theorem}

Let us indicate briefly how the map $\phi$ is defined (also
compare \cite{AR}). A $\GG$-equivariant vector superbundle $E
=E_0|E_1$, when restricted to $X^g$, becomes a spin
$\widetilde{\langle g\rangle}$-vector superbundle, where
$\widetilde{\langle g\rangle}$ denotes the subgroup of $\GG$
covering the  cyclic subgroup $\langle g\rangle$ of $G$ generated
by $g$. Using Proposition~\ref{prop_triv}, we can obtain the
following isomorphisms of $Z_G(g)$-modules:
 $$
K^-_{\widetilde{\langle g\rangle}}(X^g) \cong K(X^g) \bigotimes
R^-\left(\widetilde{\langle g\rangle}\right) \cong K(X^g)
\bigotimes R({\langle g\rangle}) \bigotimes \varepsilon_g.
 $$
In this way we obtain a map $K_{\GG}(X) \rightarrow K(X^g)
\bigotimes R({\langle g\rangle}) \otimes \varepsilon_g$. Composing
this map with the character evaluation at $g$ gives us a map
\begin{eqnarray} \label{eq_summand}
K^-_{\GG}(X) \rightarrow K(X^g) \bigotimes \varepsilon_g.
\end{eqnarray}
whose image indeed is $Z_G(g)$-invariant.

Now we claim that this map is zero when $g$ is odd (thanks to the
$\Z_2$-grading!) and thus the summation above does not involve the
odd conjugacy classes. Let $ g \in G $ be an odd element. Take an
eigenvector $v_0 + v_1$ of $g$ in the fiber of the ungraded
subbundle $E^{\mu} \subset E_0|E_1$, where $v_i \in E_i (i =
0,1)$. We see that ${g}.v_0 =\mu v_1,$ and $ {g}.v_1 =\mu v_0.$ It
follows that $v_0 -v_1$ is an eigenvector of ${g}$ with eigenvalue
$-\mu$, i.e. $v_0 -v_1 \in E^{-\mu}$. Denote by $\sharp$ the
isomorphism of $E$ which is the identity map when restricted to
$E_0$ and negative the identity map when restricted to $E_1$.
Clearly $\sharp$ sends $E^{\mu}$ to $E^{-\mu}$ and vice versa.
Thus the map (\ref{eq_summand}) becomes zero since $\mu [E^{\mu}]
+ (-\mu)[E^{-\mu}] =0$. Putting (\ref{eq_summand}) together for
all even conjugacy classes, we obtain the map $\phi$.

The rest of the proof of the above theorem is the same as in
\cite{AR} which in turn follows closely the classical case (cf.
\cite{BC, Ku, AS}).

Below we will single out a certain class of $\GG$-space $X$ with a
favorable property.
\begin{definition}
Assume we are in the setup of Theorem~\ref{th_supertech}. We say
the $\GG$-space $X$ satisfies {\em a strong vanishing property} if
for every {\em even non-split} conjugacy class $[g]$ in $G$, there
exists some element $b$ in $Z_G(g)$ such that the character
$\varepsilon_g$ of $Z_G(g)$ takes non-trivial value (which has to
be -1) and the element $b$ fixes $X^g$ pointwise.
\end{definition}

In view of Lemma~\ref{lem_chara}, if the $\GG$-space $X$ satisfies
the strong vanishing property, the isomorphism (\ref{eq_decompo})
will simplify to the following isomorphism

\begin{eqnarray*}
   \phi : \underline{K}^-_{\GG}(X)
   \stackrel{\cong}{\longrightarrow}
    \bigoplus_{[g]}
     \underline{K}(X^{g})^{Z_{G} (g)},
\end{eqnarray*}
where the summation runs over all {\em even split} conjugacy
classes in $G$.

When $X$ is a point the isomorphism $\phi$ becomes the map from a
spin supermodule of $G$ to its character. As is known \cite{Jo1,
Jo2}, the character of a spin supermodule vanishes on odd
conjugacy classes as well on even non-split conjugacy classes. In
our terminology, the one-point space for any $\GG$ automatically
satisfies the strong vanishing property. We shall see that the
strong vanishing property holds for other non-trivial examples.

\begin{remark} \rm  \label{rem_isom}
Theorem~\ref{th_supertech} contains Theorem~\ref{th_tech} as a
special case. Indeed, let $\GG = \tG_1 =\NP_1 \times \G$ for some
finite group $\G$ and let $X$ be a $\G$-space. We have an
isomorphism $$K_\G(X) \cong K^-_{\tG_1}(X), \quad V \mapsto V|V.
$$ On the other hand, we note that $\NP_1$ is isomorphic to $\Z_4
=\{1, a, z=a^2, a^3 \}$ with the $\Z_2$-grading given by letting
the generator $a$ be of degree $1$. The quotient of $\GG$ by $\{1,
z\}$ is $G = \G \times \Z_2$. The even conjugacy class in $G$ is
given by the conjugacy classes in $\G \times \{1\} = \G$.
Therefore the right-hand side in Theorem~\ref{th_supertech}
reduces to the right-hand side of Theorem~\ref{th_tech}.

It is possible to further generalize Theorem~\ref{th_supertech}
along the line of \cite{AR}.
\end{remark}
\section{The group $\Hy_n$ and $K$-theory operations}  \label{sec_oper}
In this section, we construct various $K$-theory operations based
on the finite supergroup $\Hy_n$. This is an analog of Atiyah's
construction \cite{At} of $K$-theory operations by using the
symmetric groups and implicitly Schur duality. The role of Schur
duality is replaced here by the Sergeev's generalization
\cite{Ser} of the Schur duality involving $\Hy_n$.
\subsection{The group $\Hy_n$ and a ring of
symmetric functions}
Let $V$ be a complex vector space of dimension $m$. We denote by
$q(V)$ the superalgebra of linear transformations on the
superspace $V|V =V +V$ which preserve an odd automorphism
$P:V|V\rightarrow V|V$ such that $P^2 =-1$. For example, if we
take $V =\C^m$, and $P$ to be given by the $2m\times 2m$ matrix
\begin{eqnarray*}
\left[
\begin{array}{cc} 0&I\\ -I&0\\
\end{array}
\right]
\end{eqnarray*}
then $q(V)$ is the Lie superalgebra which is obtained by the
associative superalgebra $Q(m)$ (see Section~\ref{sec_supergroup})
by taking the supercommutators.

Let us now consider the natural action of $q(V)$ on $V|V$. We may
form the $n$-fold tensor product $(V|V)^{{\otimes} n}$, on which
$q(V)$ acts naturally. In addition we have an action of the finite
supergroup $\bn$: the symmetric group $S_n$ acts on
$(V|V)^{{\otimes} n}$ by permutations with appropriate signs;
$a_i$ acts on $(V|V)^{{\otimes} n}$ by means of exchanging the
parity of $i$-th copy of $V|V$ via the odd automorphism $P$ of
$V|V$. More explicitly, $a_i$ transforms the vector
$v_1\otimes\ldots v_{i-1}\otimes v_i\otimes\ldots\otimes v_n$ in
$(V|V)^{{\otimes} n}$ into
$(-1)^{p(v_1)+\ldots+p(v_{i-1})}v_1\otimes \ldots v_{i-1}\otimes
P(v_i)\otimes\ldots\otimes v_n$. According to Sergeev \cite{Ser},
the actions of $q(m)$ and $\Hy_n$ (super)commute with each other.
Furthermore, one has

\begin{eqnarray*}
(V|V)^{{\otimes} n} \cong \sum_{\lambda}
2^{-\delta(l(\la))}U^{\lambda}_{m} \otimes T^{\lambda},
\end{eqnarray*}
where $T^{\la}$ is the irreducible $\Hy_n$-supermodule associated
to a strict partition $\la$, and $U^{\lambda}_{m}
=\Hom_{\Hy_n}(T^{\lambda}, (V|V)^{{\otimes}n})$ is an irreducible
$q(m)$-module. The expression $2^{-\delta(l(\la))}U^{\lambda}_{m}
\otimes T^{\lambda}$ above has the following meaning. Suppose that
$A$ and $B$ are two supergroups or two superalgebras and suppose
that $V_A$ and $V_B$ are irreducible supermodules over $A$ and $B$
of type $Q$, namely, ${\rm Hom}_{A}(V_A,V_A)$ and ${\rm
Hom}_B(V_B,V_B)$ are both isomorphic to the Clifford superalgebra
in one odd variable. It is known (cf. e.g. \cite{Jo1}) that
$V_A\otimes V_B$ as a module over $A\otimes B$ is not irreducible,
but decomposes into a direct sum of two isomorphic copies (via an
odd isomorphism) of the same irreducible supermodule. In our
particular setting when $l(\la)$ is odd both $T^\la$ and $U_m^\la$
are such modules (cf. \cite{Ser}). So in this case we mean to take
one copy inside their tensor product.

We introduce (cf. \cite{M}) the symmetric functions $q_n$ in the
variables $x_1, x_2, \ldots$ by the formula
\[
\sum_{n \geq 0} q_n t^n =\prod_i \frac{1+x_i t}{1-x_it}.
\]
Denote by $\Omega$ the subring of symmetric functions generated by
$q_1, q_2, q_3, \cdots,$ and $\Omega^n$ the subspace spanned by
symmetric functions in $\Omega$ of degree $n$. Put $\Omega_{\C} =
\Omega \otimes_{\Z} \C$ and $\Omega^n_{\C} = \Omega^n \otimes_{\Z}
\C$. Recall \cite{M} that a linear basis of $\Omega_{\C}$ is given
a distinguished class of symmetric functions $Q_{\la}$, call the
Schur $Q$-functions, parameterized by strict partitions $\la$.
Furthermore $\Omega^\C = \C [q_1, q_3, q_5, \cdots ]$ where $q_r$
($r$ odd) are algebraically independent. We take the convention
that when the set of variables is finite, say $x =(x_1, x_2,
\ldots, x_m)$, we set $Q_{\la}(x) = Q_{\la}(x_1, x_2, \ldots, x_m,
0, 0 \ldots)$. According to Sergeev \cite{Ser}, the trace of the
diagonal matrix $D={\rm diag} (x_1,\ldots,x_m; x_1,\ldots,x_m)$ in
$q(m)$ acting on ${U^\la_{m}}$ is
\begin{eqnarray} \label{eq_trace}
Tr D|U^\la_{m} = 2^{\frac{\delta(l(\la))-l(\la)}{2}}Q_\la(x).
\end{eqnarray}

Given a spin $\bn$-supermodule $W$, we define $W(V)$ to be the
space of $\Hy_n$-invariants $ \left( W\bigotimes (V|V)^{\otimes n}
\right)^{\bn}. $ It is easy to check that the correspondence $V
\mapsto W(V)$ is functorial on $V$. In particular, if we take a
diagonalizable linear transformation $l: V \rightarrow V$ with
eigenvalues $x_1, \ldots, x_m$, then the eigenvalues of the
induced map $W(l) = Id_W \otimes (l|l)^{\otimes n}: W(V)
\rightarrow W(V)$ are monomials in $x_1, \ldots, x_m$ of degree
$n$. In particular the trace of $ W(l)$ is a symmetric polynomial
in $x_1, \ldots, x_m$ of degree $n$ with integer coefficients. One
can argue that this symmetric polynomial (for $m \geq n$) is the
restriction of a unique symmetric function in infinite many
variables.

By definition we have $(W_1 \bigoplus W_2)(V) = W_1(V) \bigoplus
W_2(V).$ It follows the mapping $W \mapsto Tr W(l)$ induces a map,
denoted by $ch$, from $R^-(\bn)$ to the space of symmetric
functions of degree $n$. Note that $W(V) \cong \Hom_{\Hy_n} (W,
(V|V)^{\otimes n})$ since all the character values of $W$ are
real, i.e. $W$ is self-dual. It follows from (\ref{eq_trace}) that
if we take $W$ to be $T^{\la}$ then $ch$ sends the class function
associated with $T^{\la}$ to
$2^{\frac{\delta(l(\la))-l(\la)}{2}}Q_\la(x)$. In this way we have
defined a map $ch$ from $R^- = \oplus_n R^-(\Hy_n)$ to the ring
$\Lambda_{\C}$ of symmetric functions.

\begin{remark} \rm
It is possible that one can develop this approach coherently to
study the map $ch: R^- \rightarrow \Lambda_\C$ without referring
to the Lie superalgebra $q(V)$ and thus essentially independent of
the work of Sergeev \cite{Ser}, as sketched below. This is a super
analog of an approach adopted in Knutson \cite{Kn} in the setup of
symmetric groups. For example, we can start by arguing that the
space $W(V)$ associated to the basic spin supermodule $L_n$ (with
character $\xi_n$) is the $n$-th supersymmetric algebra
$\sss^n(V|V) = \oplus_{i=0}^n S^i V \otimes \Lambda^{n-i}V$, and
thus $ch (\xi^n) = \sum_{i=0}^n h_i e_{n-i} =q_n$ (cf. Ex. 6(c),
pp261, \cite{M}), which is exactly the Schur $Q$-function
$Q_{(n)}$. Here $h_r$ and $e_r$, which stand for the {\it $r$-th
elementary} and {\it complete} symmetric function respectively,
are the traces of ${\rm diag} (x_1,\ldots,x_m)$ on $S^r V$ and
$\Lambda^r V$. One can easily show by using the Frobenius
reciprocity that $ch:R^-(\bn) \rightarrow \Lambda$ is an algebra
homomorphism. The image of $ch$ contains $\Omega^n_\C$ as we have
just seen that it contains $q_n$ for all $n$. By comparing the
dimensions, the characteristic map $ch: R^- \rightarrow \Omega_\C$
is indeed an isomorphism.
\end{remark}

There is another way to define the characteristic map as follows
(cf. \cite{Jo2, Ser}). Denote by $p_n$ the $n$-th power sum
symmetric functions, and define $p_{\mu} =p_{\mu_1}p_{\mu_2}
\ldots$ for a partition $\mu =(\mu_1, \mu_2, \ldots)$. Given a
character $\chi \in R^-$, we define the characteristic map $ch' :
R^- \rightarrow \Omega_{\C}$ by
\begin{eqnarray*}
ch'(\chi) =\sum_{\mu \in {\cal OP}_n} z_{\mu}^{-1}
\chi_{\mu}p_{\mu},
\end{eqnarray*}
where $ z_{\mu} = \prod_{i\geq 1}i^{m_i}m_i!$ for $\mu =(1^{m_1}
3^{m_3} \ldots )$ and $\chi_{\mu}$ the character value at the even
split conjugacy class $D^+_{\mu, \emptyset}$. It is known
\cite{Jo2} that $ch'$ is an algebra isomorphism and it sends
$\xi^n$ to $q_n$ for all $n$. Thus $ch'$ coincides with the
characteristic map $ch$ we defined above.
\subsection{K-theory operations}
Let $\GG$ is a finite supergroup with $G = \GG/\langle 1,z
\rangle$, and let $X$ be a trivial $G$-space (and thus
$\GG$-space). Given a $G$-supermodule $V$ in $R_{\Z}^-(G)$ and a
complex vector bundle $E$ on $X$, $E \otimes V$ can be given a
natural $\GG$-equivariant vector superbundle structure over $X$ by
letting $\GG$ act on only the factor $V$. Obviously $E \otimes V$
lies in the category ${\cal C}_{\GG}(X)$. In this way we obtain a
map $R_{\Z}^-(\GG) \otimes K(X) \longrightarrow K^-_{\GG}(X)$.
Then we have the following (compare \cite{Seg1}).
\begin{proposition} \label{prop_triv}
Under the above setup, there is a canonical isomorphism
 $$ K^-_{\GG}(X)\stackrel{\cong}{\longrightarrow} K(X) \bigotimes R_{\Z}^-(\GG).
 $$
\end{proposition}

One needs some extra care to define the inverse of the isomorphism
in the above proposition. Given a $G$-supermodule $V$ in
$R_{\Z}^-(G)$ and a complex vector bundle $E$ on $X$,  consider
$\Hom_{G} ( {\bf V}, E)$, where $\bf V$ is the trivial
$G$-superbundle over $X$ associated to $V$. $\Hom_{G} ( {\bf V},
E)$ is isomorphic to $E$ if $V$ is of type $M$, but isomorphic to
$E \otimes C_1 =E|E$ if $V$ is of type $Q$, where $C_1$ is the
Clifford algebra in one variable. It is perhaps more natural to
replace $K(X)$ in the proposition above by an isomorphic space
$K^-_{\NP_1}(X)$, cf. Remark~\ref{rem_isom}.

Below we will construct various $K$-theory operations based on the
construction in the previous subsection. It is a super analog of
an approach due to Atiyah \cite{At} who used the symmetric group
representations.

Let $E$ be a vector bundle over $X$, consider the $n$-th tensor
power $(E|E)^{\otimes n}$ of the vector superbundle $E|E$. The odd
operator $P$ ($P^2= -1$) acts on each factor $E|E$ fiberwise and
induces an action of the finite supergroup $\NP_n$ on $(E
|E)^{\otimes n}$. The symmetric group $S_n$ also acts on $(E
|E)^{\otimes n}$ in a natural way. The joint action of $\NP_n$ and
$S_n$ then gives rise to an action of the finite supergroup $\bn$
on $(E |E)^{\otimes n}$. We have the following decomposition
\begin{eqnarray*}
 (E |E)^{\otimes n} \cong
 \sum_{\pi}2^{- \delta(l(\pi))} T^{\pi} \bigotimes \pi (E),
\end{eqnarray*}
where $\pi$ is a strict partition of $n$, and $\pi (E)
=\left(T^{\pi}\otimes (E|E)^{\otimes n} \right)^{\bn}$ is a vector
bundle on $X$. Clearly one can extend the definition of the vector
bundle $\pi (E)$ associated to any spin supermodule $\pi$ of $\bn$
so that it is additive on $\pi$. In this way we obtain a ring
homomorphism $j^-: R^- \rightarrow {\rm Op}(K)$, where $R^-
=\bigoplus_n R^-(\widetilde{H}_n)$ and ${\rm Op}(K)$ is the ring
of $K$-theory operations on $K(X)$ (cf. \cite{At}).

We note that if $\pi$ is the one-part partition $(n)$ then $\pi
(E)$ is the $n$th supersymmetric power $\sss^n$:
\begin{eqnarray} \label{eq_susy}
\sss^n (E|E) =\sum_{i=0}^n S^i E \otimes \Lambda^{n-i}E.
\end{eqnarray}
For odd $n$, we denote by $\psi^n$ the operation corresponding to
the class function $\sigma_n \in R^-(\bn)$ which takes value $n$
in the even split conjugacy class of type $((n), \emptyset)$ and
zero elsewhere. Denote by $\sss_t(E) = \sum_{n =0}^{\infty} \sss^n
(E|E) t^n,$ where $t$ is a formal variable.

\begin{proposition}
The operation $\psi^r$ ($r$ odd) coincides with the usual $r$th
Adams operation. In particular it is additive. Furthermore, we
have
\begin{eqnarray}  \label{eq_exp}
\sss_t(E) = \exp \left( \sum_{r >0, \; odd} 2 \psi^r(E) t^r/{r}
\right).
\end{eqnarray}
\end{proposition}

\begin{demo}{Proof}
We first prove (\ref{eq_exp}). Since the operations $\sss^n$'s and
$\psi^r$'s are obtained from the ring homomorphism $j^-: R^-
\rightarrow {\rm Op}(K)$, we only need to show that the
corresponding identity holds in $R^-$, or alternatively in
$\Omega_\C$. But this is a classical identity \cite{M}
\[
\sum_{n \geq 0} q_n t^n =\exp \left( \sum_{r >0, \; odd} 2 p_r
t^r/{r} \right).
\]

Now we denote by $\sigma_t =\sum_n \sigma^n t^n, \Lambda_t =\sum_n
\Lambda^n t^n$. Let us denote by $\tilde{\psi}^r$ the $r$th Adams
operations for the time being. It is classical that
\begin{eqnarray*}
\sigma_t(E) &=&\exp \left( \sum_{r >0} \tilde{\psi}^r(E) t^r/{r}
\right) \\
 \Lambda_t(E) &=&\exp \left(\sum_{r >0}(-1)^{r-1}
\tilde{\psi}^r(E) t^r/{r} \right).
\end{eqnarray*}
We have from (\ref{eq_susy}) that $\sss_t(E) = \Lambda_t(E)
\sigma_t(E)$ and therefore $$\sss_t(E) =\exp \left( \sum_{r >0,
\;odd} 2\tilde{\psi}^r(E) t^r/{r} \right).$$ Comparing with
(\ref{eq_exp}) which we have already established, we have $\psi^r
=\tilde{\psi}^r $.
\end{demo}

It follows that $Q_t(E+F) =Q_t(E)Q_t(F)$ and $Q_t(E-F)
=Q_t(E)Q_{-t}(F),$ where $ E,F \in K(X)$, since the Admas
operations are additive. For example, the second equation reads
componentwise as follows:
\begin{eqnarray*}
 Q^n(E-F|E-F) &=& \bigoplus_{i=0}^n (-1)^i
Ind^{\Hy_n}_{\Hy_{n -i} \widehat{\times} \Hy_i} (Q^{n -i}(E|E)
\otimes Q^i(F|F)).
\end{eqnarray*}

\begin{remark} \rm
The ring $\Lambda$ of symmetric functions is a basic model for
(free) $\la$ rings, and indeed $\la$ rings can be defined by
axiomizing various properties of natural operations on $\Lambda$,
where the Adams operations play a crucial role (cf. \cite{Kn}). We
can define similarly a notion of $Q$-$\la$ ring with Adams
operations of odd degrees only, using $\Omega_\C$ as a basic
model. Then what we have just shown is that $K(X)$ admits a
$Q$-$\la$ ring structure. We will see the structure of a $Q$-$\la$
ring instead of a $\la$ ring shows up naturally in some fairly
non-trivial setup in Section~\ref{sec_structure}.
\end{remark}
\section{Algebraic structures on $\tFG$} \label{sec_structure}

In this section, we shall study in detail the $K$-groups
$\kbgnfree$ of generalized symmetric products for all $n$
simultaneously.
\subsection{Generalized symmetric products}
Our main examples in this paper are as follows. Let $\G$ be a
finite group and let $X$ be a $\G$-space. The $n$-th Cartesian
product $X^n$ is acted by the finite supergroup $\tG_n = \NP_n
\rtimes \Gn$ in a canonical way: $\NP_n$ acts trivially on $X^n$;
$\G^n$ acts on $X^n$ factorwise while $S_n$ by permutations, and
this gives rise to a natural action of the wreath product $\Gn
=\G^n \rtimes S_n$ by letting
\begin{eqnarray*}
 a . (x_1, \ldots, x_n)
  = (g_1 x_{s^{-1} (1)}, \ldots, g_n x_{s^{-1} (n)})
\end{eqnarray*}
where $ a = ( (g_1, \ldots, g_n), s) \in \Gn$, and $x_1, \ldots,
x_n \in X$. Note that orbifolds $X^n/S_n$ are often called
symmetric products. We may refer to $X^n/\Gn$, or $X^n$ with
$\Gn$-action, or rather $X^n$ with $\tG_n$-action as {\em
generalized symmetric products}.

Our earlier general construction when applied to the generalized
symmetric products gives us the category ${\cal C}^-_{\tG_n}
(X^n)$ and its associated $K$-group $\underline{K}^-_{\tG_n}
(X^n)$. It turns out this category affords an equivalent
description below which affords more transparent geometric
meaning.

The wreath product $\Gn$ acts on the vector space $\C^n$ naturally
by letting $\G$ act trivially and $S_n$ act as the permutation
representation. This action preserves the standard quadratic form
on $\C^n$. We denote by ${\bf n}$ such a $\Gn$-vector bundle $X^n
\times \C^n$ over $X^n$. We denote by $C_n$ the complex Clifford
algebra associated to $\C^n$ and the standard quadratic form on
it. The action of $\Gn$ on $\C^n$ induces a natural action on
$C_n$. We denote by $C({\bf n})$ the associated $\Gn$-vector
(super)bundle $X^n \times C_n$ on $X^n$, which is the Clifford
module on $X^n$ associated to the vector bundle ${\bf n}$.

We introduce the following category ${\cal C}_{\Gn}^{\bf n}(X^n)$:
the objects consist of complex vector superbundles $E = E_0 + E_1$
on $X^n$ equipped with compatible actions of $\Gn$ and the
Clifford algebra $C_n$ associated to $\C^n$ with the standard
quadratic form. That is, $E$ is a $\Z_2$-graded $C({\bf
n})$-module and a $G$-equivariant vector bundle over $X^n$ such
that
\begin{eqnarray}  \label{eq_compat}
g.(v. \xi) =(g.v).(g.\xi), \quad g\in\Gn, v \in {\bf n},\xi\in E.
\end{eqnarray}
Given two superbundles $E, F$ in ${\cal C}_{G}^{\bf n}(X^n)$, the
space of $(G, C({\bf n}))$-equivariant homomorphisms of vector
superbundles admits a natural $\Z_2$-gradation. Since the twisted
group algebra of $\NP_n$ is isomorphic to the Clifford algebra
$C_n$, the category ${\cal C}_{G}^{\bf n}(Y)$ is then obviously
equivalent to the category ${\cal C}^-_{\tG_n}(X^n)$, and so are
the corresponding $K$-groups.

When $X$ is a point, ${K}^-_{\tG_n}(pt) = {K}^{-,0}_{\tG_n}(pt)$
reduces to the Grothendieck group $R^-_{\Z}(\tG_n)$ of the spin
supermodules of $\tG_n$.

\begin{remark}  \rm
We may replace the rank $n$ vector bundle $\bf n$ over $X^n$ above
by the $n$th direct sum of a non-trivial line bundle endowed with
a quadratic form, and modify the construction of the category
${\cal C}_{\Gn}^{\bf n}(X^n)$ accordingly. We conjecture that the
resulting $K$-group is isomorphic to ${K}^-_{\tG_n}(X^n)$.
\end{remark}

\begin{remark} \rm
We may reverse the above consideration in a more general setup as
below. Take a $G$-space $Y$ and a $G$-superbundle $C({\bf n})$
over $Y$ whose fiber is the Clifford algebra in $n$ variables.
Assume that there exists $n$ sections $a_1, \ldots, a_n$ of the
superbundle $C({\bf n})$ which fiberwise generate the Clifford
algebra and $G$ permutes $a_1, \ldots, a_n$. It is of interest in
algebraic topology (cf. Karoubi \cite{Ka}) to study the category
${\cal C}_{G}^{\bf n}(Y)$ of $G$-equivariant vector superbundles
over $Y$ which are compatible with the the Clifford module
structure in the sense of (\ref{eq_compat}) and study its
associated $K$-group. Then we may form the finite supergroup
$\NP_n\rtimes G$ and reinterpret the category ${\cal C}_{G}^{\bf
n}(Y)$ as the category ${\cal C}^-_{\NP_n\rtimes G}(Y)$, and then
apply Theorem~\ref{th_supertech} to the study of the associated
$K$-group. In particular if we take $Y=X^n$ and $G=\Gn$ for a
$\G$-space $X$, we recover our main examples of generalized
symmetric products.
\end{remark}
\subsection{Hopf algebra $\tFG$}
Let $\widetilde{H}$ be a $\Z_2$-graded subgroup of  a finite
supergroup $\GG$ with the same distinguished even element $z$, and
let $X$ be a $G$-space $X$ which is regarded as a $\GG$-space,
where $G =\GG/\langle 1,z \rangle $. we can define the restriction
map $Res^G_H: \underline{K}^-_{G} (X) \longrightarrow
\underline{K}^-_{H} (X)$ and the induction map $Ind^G_H:
\underline{K}^-_{H} (X) \longrightarrow \underline{K}^-_{G} (X)$
in the same way as in the usual equivariant $K$-theory. When it is
clear from the text, we will often abbreviate $Res^G_H$ as $Res_H$
or $Res$. Similar remarks apply to the induction map.

Now we introduce the direct sum of equivariant $K$-groups
\begin{eqnarray*}
\tFG = \bigoplus_{n \geq 0} \kbgnfree, \quad
 \tFGvar =\bigoplus_{n \geq 0} t^n \kbgnfree,
\end{eqnarray*}
where $\tG_0$ is the one-element group by convention, and $t$ is a
formal variable counting the graded structure of $\tFG$. We also
set
 $$ \dim_t \tFG =\sum_{n \geq 0} t^n \dim
\kbgnfree. $$

We define a multiplication $\cdot$ on  the space $ \tFG$ by a
composition of the induction map and the K\"unneth isomorphism
$k$:
\begin{eqnarray}  \label{eq_mult}
   \underline{K}^-_{\tG_n} (X^n) \bigotimes  \underline{K}^-_{\tG_m} (X^m)
      \stackrel{k }{\longrightarrow}
       \underline{K}^-_{\tG_n \widehat{\times} \tG_m} (X^{n +m} )
 \stackrel{\mbox{Ind}}{\longrightarrow}
  \underline{K}^-_{\tG_{n+m}} (X^{n +m} ).
\end{eqnarray}
We denote by $ 1 $ the unit in $\underline{K}^-_{\tG_0} (X^0)$
which can be identified with $ \C $.

On the other hand we can define a comultiplication $\Delta$ on
$\tFG$ to be a composition of the inverse of the K\"unneth
isomorphism and the restriction from $\tG_n$ to $ \tG_m
\widehat{\times} \tG_{n -m}$:
\begin{eqnarray} \label{eq_comult}
   \underline{K}^-_{\tG_n} (X^n)
 \longrightarrow \bigoplus_{m =0}^n
        \underline{K}^-_{\tG_{m } \widehat{\times} \tG_{n -m} } (X^n)
 \stackrel{k^{-1}}{ \longrightarrow} \bigoplus_{m =0}^n
    \underline{K }^-_{\tG_{m }} (X^m) \otimes
    \underline{K}^-_{\tG_{n -m}} (X^{n -m}).
\end{eqnarray}
We define the counit $\epsilon : \tFG \longrightarrow \C$ by
sending $\underline{K}^-_{\tG_n} (X^n)$ $(n >0)$ to $0$ and $ 1
\in K^-_{\tG_0} (X^0) = \C $ to $1$.

\begin{theorem} \label{th_hopf}
With various operations defined as above, $ \tFG$ is a graded Hopf
algebra.
\end{theorem}

The proof is the same as the proof of the Hopf algebra structure
on a direct sum of equivariant $K$-groups $\oplus_n
\underline{K}_{\G_n}(X^n)$ \cite{W1}, where a straightforward
generalization to the equivariant $K$-groups of the Mackey's
theorem plays a key role. One easily checks the super version of
the Mackey's theorem can be carried over to our $K$-group setup.
In the case when $X$ is a point and thus $\underline{K}^-_{\tG_n}
(X^n) \cong R^-(\tG_n)$, the Hopf algebra structure above is also
treated in \cite{JW}.
\subsection{Description of the algebra $\tFG$} \label{subsect_alg}
Take an even split element $a =(\g, \sigma)$ in $\bG_n$ of type
$\rho =(\rho^+, \emptyset)$, where $g=(g_1, \cdots, g_n)$, $g_i =
(\alpha_i, \varepsilon_i) \in\G\times \Z_2$, and $\rho^+
=(m_r^+(c))_{r,c} \in {\cal OP}_n(\G_*)$. We define $\overline{a}
=((\alpha_1, \cdots, \alpha_n), \sigma) \in \Gn = \G^n \rtimes
S_n$. Since the subgroup $\NP_n \leq \bG_n$ acts on $X^n$
trivially, the orbit space $(X^n)^a /Z_{\bG_n}(a) $ is identified
with $ (X^n)^{\overline{a}} /Z_{\Gn}(\overline{a})$ which has been
calculated earlier in Lemmas 4 and 5 of \cite{W1}. We make a
convention here to denote the centralizer $Z_G (g)$ (resp. $X^g$,
$X^g /Z_G (g)$) by $Z_G (c)$ (resp. $X^c$, $X^c /Z_G (c)$) by
abuse of notations when the choice of a representative $g$ in a
conjugacy class $c$ of a group $G$ is irrelevant. For a fixed $c
\in \G_*$, recall that $c_n$ ($n \in 2\Z_+ +1$) is the even split
conjugacy class in $\bG_n$ of the type $(\rho^+, \emptyset)$,
where the partition-valued function $\rho^+$ takes value the
one-part partition $(n)$ at $c \in \G_*$ and zero elsewhere.

\begin{lemma} \label{lem_orbit}
Let $a \in \bG_n$ be an even split element of type $\rho= (\rho^+,
\emptyset)$, where $\rho^+ =(m_r^+(c))_{r\geq 1,c \in \G_*}
\in{\cal OP}_n(\G_*)$. Then the orbit space $( X^n )^a
/Z_{\bG_n}(a)$ can be naturally identified with
 $$ \prod_{r,c} S^{m^+_r (c)} \left( X^{c} / Z_\G (c) \right), $$
where $S^{m}(\cdot)$ denotes the $m$-th symmetric product. In
particular, the orbit space $( X^n )^{c_n} /Z_{\bG_n}(c_n)$ can be
naturally identified with $ X^{c}/Z_\G (c).$
\end{lemma}

In view of Lemma~\ref{lem_chara}, let us recall how we have
classified the split conjugacy classes of $\bG_n$ in Theorem 1.2
of \cite{JW}. In order to show that a given element $a$ in $\bG_n$
is non-split, an explicit element, say $b$, is constructed in the
centralizer $Z_{\bG_n}(a)$ such that the character $\varepsilon_a$
takes value $-1$ at $b$. This was achieved in \cite{JW} case by
case. On the other hand, we observe that in all cases the element
$b$ fixes $(X^n)^a$ pointwise! In other words, we have the
following.

\begin{lemma} \label{lem_vanishing}
The $\tG_n$-space $X^n$ satisfies the strong vanishing property.
\end{lemma}

We now give an explicit description of $\tFG =\bigoplus_{n \geq 0}
\kbgnfree$ as a graded algebra.
\begin{theorem} \label{th_size}
As a $(\Z_+ \times \Z_2)$-graded algebra ${\cal F}^-_{\G}(X, t)$
is isomorphic to the supersymmetric algebra
  $ {\cal S} \left( \bigoplus_{ r=1}^{\infty} t^{2r-1} \underline{K}_\G(X)
      \right)$. In particular, we have
 $$
 \dim_t \tFG =
 \frac{ \prod_{ r=1}^{\infty}
 (1 + t^{2r-1})^{ \dim K^1_\G (X)} }{ \prod_{ r=1}^{\infty}
    (1 - t^{2r-1})^{ \dim K^0_\G (X)} }.
 $$
\end{theorem}
Here the supersymmetric algebra is equal to the tensor product of
the symmetric algebra $ S \left( \bigoplus_{ r=1}^{\infty}
t^{2r-1} \underline{K}^0_\G(X)
      \right)$ and the exterior algebra
$\Lambda \left( \bigoplus_{ r=1}^{\infty} t^{2r-1}
\underline{K}^1_\G(X) \right)$.

\begin{demo}{Proof}
Take an even split element $a  \in \bG_n$ of type $\rho= (\rho^+,
\emptyset)$, where $\rho^+ =(m_r^+(c))_{r\geq 1,c \in \G_*}
\in{\cal P}_n(\G_*)$. By Lemma~\ref{lem_orbit} and the K\"unneth
formula, we have
 \begin{eqnarray}
   \underline{K}( (X^n)^a / Z_{\bG_n} (a) )
   & \approx &
   \bigotimes_{c \in \G_*, r \geq 1 \; odd}
    \left(
      \left( \underline{K}(X^{c})^{ Z_\G (c)}
      \right)^{\bigotimes m_r (c)}
    \right)^{S_{m^+_r (c)}  }  \nonumber\\
   & \approx & \bigotimes_{c \in \G_*, r \geq 1 \; odd}
   {\cal S}^{m^+_r (c)} ( \underline{K}(X^{c} / Z_\G (c)) ).   \label{eq_symm}
 \end{eqnarray}

We now calculate as follows. The statement concerning $\dim_t
{\tFG} $ follows from this immediately.

 \begin{eqnarray*}
 {\cal F}^-_\G (X, t)
   & \cong &   \bigoplus_{n \geq 0} t^n
     \bigoplus_{[a] \in (\bG_n)_*^{e.s.}}
      \underline{K} \left( (X^n)^a / Z_{\bG_n} (a) \right)
       \\
        & & \qquad\qquad\qquad\qquad\qquad\qquad\qquad
        \mbox{by Theorem } \ref{th_supertech} \mbox{ and Lemma}~\ref{lem_vanishing},  \\
   & \cong &
   \bigoplus_{ n \geq 0}
     \bigoplus_{ \{m^+_r(c)\}_{c,r} \in {\cal OP}_n ( \G_*)}
        t^n \bigotimes_{c, r \; odd}
        {\cal S}^{m^+_r (c)} (\underline{K} (X^{c} / Z_\G (c)) ) \\
        & & \qquad\qquad\qquad\qquad\qquad\qquad\qquad\qquad\quad
        \mbox{by }(\ref{eq_symm})
        \mbox{ and Theorem } \ref{th_splitclass}, \\
   & \cong & \bigoplus_{\{m^+_r(c)\}_{c,r} \in {\cal OP} ( \G_*)}
        \bigotimes_{c, r\; odd} {\cal S}^{m^+_r (c)}
          \left(t^r \underline{K}(X^{c} / Z_\G (c)) \right)
            \nonumber     \\
   & \cong & \bigoplus_{\{m_r\}_{r} } \bigotimes_{ r \geq 1\; odd}
   {\cal S}^{m_r}
          \left( \bigoplus_{c \in \G_*} t^r
          \underline{K} (X^{c} / Z_\G (c))
          \right)
       \quad\mbox{ where $m_r = \sum_{c} m^+_r (c)$},
                               \label{eq_sum}     \\
   & \cong & \bigoplus_{\{m_r\}_{r} } \bigotimes_{ r \geq 1 \; odd}
             {\cal S}^{m_r} (t^r \underline{K}_\G (X) )
    \; \quad\quad\quad\quad\quad\quad\quad\quad
    \quad\;\mbox{by Theorem } \ref{th_tech}, \label{eq_again}       \\
   & \cong & {\cal S} \left( \bigoplus_{ r=1}^{\infty}
t^{2r-1}
    \underline{K}_\G(X) \right).
                  \nonumber
 \end{eqnarray*}
\end{demo}

Recall that the orbifold Euler number $e(X,\G)$ was introduced by
Dixon, Harvey, Vafa and Witten \cite{DHVW} in the study of
orbifold string theory. It is subsequently interpreted as the
Euler number of the equivariant $K$-group $K_\G(X)$, cf. e.g.
\cite{AS}. If we define the Euler number of the generalized
symmetric product to be the difference
 $$
 e(X^n,\tG_n) := \dim K^{-,0}_{\tG_n}(X^n) - \dim
 K^{-,1}_{\tG_n}(X^n),
 $$
then we obtain the following corollary.

\begin{corollary}  \label{cor_euler}
The Euler number $e(X^n,\tG_n)$ is given by the following
generating function:
 $$
\sum_{n=0}^{\infty} e(X^n,\tG_n) t^n =\prod_{ r=1}^{\infty} (1
-t^{2r-1})^{ -e(X,\G)}.
 $$
\end{corollary}

In the case when $X$ is a point, we obtain the following corollary
(also cf. \cite{JW}).

\begin{corollary}
When $X$ is a point and thus $\tFG \cong R^-_{\G}$, we have
\begin{eqnarray*}
\sum_{n \geq 0} t^n \dim R^-(\tG_n)
 &=& \prod_{ r=1}^{\infty} (1 -t^{2r-1})^{ -|\G_*|} .
\end{eqnarray*}
\end{corollary}
\subsection{Twisted vertex operators and $\tFG$} \label{subsect_twistvo}
In the following, we will define various
$K$-theory maps appearing in the following diagram ($n$ odd):
 \begin{eqnarray} \label{eq_maps}
  \underline{K}_\G (X) & \stackrel{\aleph}{\longrightarrow} &
   \underline{K}^-_{\tG_1} (X)
     \stackrel{{\boxtimes}n}{\longrightarrow} \underline{K}^-_{\tG_n}(X^n)
     \nonumber \\
     & \stackrel{\phi_n}{\longrightarrow}&
    \bigoplus_{[a] \in (\bG_n)_*^{e.s}}
      \underline{K} \left( (X^n)^a / Z_{\bG_n} (a) \right) \nonumber \\
     & \stackrel{\stackrel{pr}{\rightleftarrows} }{\iota}&
    \bigoplus_{c \in \G_*}
      \underline{K} \left( (X^n)^{c_n} / Z_{\bG_n} (c_n) \right) \\
    & \stackrel{\vartheta}{\longrightarrow} &
    \bigoplus_{c \in \G_*} \underline{K} \left( X^{c} / Z_\G (c)
    \right)
    \stackrel{\phi}{\longleftarrow} \underline{K}_\G (X).\nonumber
 \end{eqnarray}

Noting that $\tG_1 = \G \times \Pi_1$, we have a canonical
isomorphism, denoted by $\aleph$, from $K_\G (X)$ to
$K_{\tG_1}^-(X)$ given by $E \mapsto E|E$. Given a
$\G$-equivariant vector bundle $V$, consider the $n$-th outer
tensor product $(V|V)^{\boxtimes n}$ which is a vector superbundle
over $X^n$. The odd operator $P$ acting on each factor $V|V$
induces an action of the finite supergroup $\NP_n$ on
$(V|V)^{\boxtimes n}$ while the wreath product $\Gn$ acts on
$(V|V)^{\boxtimes n}$ by letting
\begin{eqnarray}  \label{eq_act}
 ( (g_1, \ldots, g_n), s). u_1 \otimes \ldots \otimes u_n
 = g_1 u_{ s^{-1}(1)} \otimes \ldots \otimes  g_n u_{ s^{-1}(n)}
\end{eqnarray}
where $g_1, \ldots, g_n \in \G, s \in S_n $ and $u_i \in V|V, i=1,
\ldots, n$. It is easy to check the combined action gives rise to
an action of the finite supergroup $\tG_n$ on $(V|V)^{\boxtimes
n}$, and $(V|V)^{\boxtimes n}$ endowed with such an $\tG_n$ action
is an $\tG_n$-equivariant vector superbundle over $X^n$. On the
other hand we can define an $\tG_n$ action on $ V^{\boxtimes n}
\otimes (\C^{1|1})^{\boxtimes n}$ as follows: $\Gn$ acts on the
first factor $V^{\boxtimes n}$ only while $\NP_n$ acts only on the
second factor $(\C^{1|1})^{\boxtimes n}$; the symmetric group
$S_n$ acts diagonally. One can check that the combined action
gives $ V^{\boxtimes n} \otimes (\C^{1|1})^{\boxtimes n}$ the
structure of an $\tG$-equivariant vector superbundle over $X^n$.
We easily see that $(V|V)^{\boxtimes n}$ is canonically isomorphic
to $V^{\boxtimes n} \otimes (\C^{1|1})^{\boxtimes n}$ as a
$\tG_n$-equivariant superbundle.

\begin{remark} \rm \label{rem_tensor}
Note that the above $(\C^{1|1})^{\boxtimes n}$ is precisely the
basic $\tG_n$-supermodule $L_n$. In general for a given
$\tG_n$-supermodule $M$ with character $\chi$, we can define an
$\tG_n$-equivariant superbundle structure on $V^{\boxtimes n}
\otimes M$ when replacing $L_n$ above by $M$. We will write the
corresponding element in $K_{\tG_n}(X^n) $ as $V^{\boxtimes n}
\star \chi$. This defines an additive map from $R(\Gn)$ to
$K^-_{\tG_n}(X^n)$ by sending $\chi$ to $V^{\boxtimes n} \star
\chi$.
\end{remark}

Sending $V|V$ to $(V|V)^{\boxtimes n}$ gives rise to the K-theory
map ${\boxtimes}n$. More explicitly, given $V, W $ two
$\G$-equivariant vector bundles on $X$, we use $V$ itself to
denote the corresponding element in $K_\G(X)$ by abuse of
notation. Then
\begin{eqnarray} \label{eq_comp}
 (V|V -W|W)^{\boxtimes n} = \sum_{j =0}^n (-1)^j
  \Ind_{\tG_{n -j} \widehat{\times} \tG_j}^{\tG_n}
  \left((V|V)^{\boxtimes (n -j)} \boxtimes (W|W)^{\boxtimes j}\right).
\end{eqnarray}
Here $(V|V)^{\boxtimes (n -j)}$ and $(W|W)^{\boxtimes j}$ carry
the standard actions of $\tG_{n-j}$ and respectively $\tG_j$. The
map $\phi_n$ is the isomorphism in Theorem~\ref{th_supertech}
given by the summation $\sum_{\rho^+ \in {\cal OP}_n(\G_*)}
(\phi_n)_{\rho^+}$ over the even split conjugacy classes of
$\bG_n$ of type $(\rho^+, \emptyset)$ when applying to the case
$X^n$ with the action of $\tG_n$. The map $pr$ is the projection
to the direct sum over the even split conjugacy classes $c_n$ of
$\bG_n$ while $\iota$ denotes the inclusion map. The map
$\vartheta$ denotes the natural identification given by
Lemma~\ref{lem_orbit}. Finally the last map $\phi$ is the
isomorphism given in Theorem~\ref{th_tech}.

We introduce in addition the following $K$-theory operations.

\begin{definition} \label{def_operat}
For $n \in 2\Z_+ +1$, we define the following $K$-theory
operations as composition maps:
\begin{eqnarray*}
  \psi^n &:= &
  n \phi^{-1} \circ \vartheta \circ pr \circ \phi_n
  \circ {\boxtimes}n \circ \aleph:
    \underline{K}_\G(X) \longrightarrow  \underline{K}_\G(X),  \\
 \varphi^n &:= &
  n \phi_n^{-1} \circ \iota \circ {pr}
   \circ \phi_n \circ {\boxtimes}n :
    \underline{K}^-_{\tG_1}(X) \longrightarrow \underline{K}^-_{\tG_n}(X^n), \\
 {ch}_n &:= &  \phi^{-1} \circ \vartheta \circ pr \circ \phi_n:
   \underline{K}^-_{\tG_n}(X^n) \longrightarrow \underline{K}_\G(X),  \\
 \varpi_n &:= & \phi_n^{-1} \circ \iota
 \circ \vartheta^{-1} \circ \phi:
    \underline{K}_\G(X) \longrightarrow \underline{K}^-_{\tG_n}(X^n).
\end{eqnarray*}
\end{definition}

Recall that the notation $\psi^n$ ($n$ odd) was used in
Section~\ref{sec_oper} to denote the $n$th Adams operation. We
shall see the $\psi^n$ defined here for  $\G$ trivial coincides
with the $n$th Adams operation tensored with $\C$. This is why we
have chosen to use the same notation. We list some properties of
these K-theory maps which follows directly from defintions.

\begin{proposition} \label{prop_property}
The following identities hold:
 \begin{eqnarray*}
  ch_n \circ \varpi_n  =I|_{\underline{K}_\G(X)},\quad
  \varpi_n \circ \psi^n  = \varphi^n  \circ \aleph,\quad
  ch_n \circ \varphi^n  \circ \aleph =\psi^n,
 \end{eqnarray*}
where $I|_{\underline{K}_\G(X)}$ denotes the identity operator on
$\underline{K}_\G(X)$.
\end{proposition}

\begin{lemma}
Both $\psi^n$ and $\varphi^n$ ($n$ odd) are additive K-theory
maps. In particular, for $\G$ trivial, the operation $\psi^n$
given in the definition~\ref{def_operat} coincides with the $n$th
Adams operations on $\underline{K}(X)$.
\end{lemma}

\begin{demo}{Proof}
We sketch a proof. By definition $\varpi_n$ is additive and
$\aleph$ is an isomorphism. Thanks to the equality $\varpi_n \circ
\psi^n = \varphi^n \circ \aleph$ by
Proposition~\ref{prop_property}, to show that $\varphi^n$ is
additive, it suffices to check that $\psi^n$ is additive. This can
be proved in a parallel way by using (\ref{eq_comp}) as Atiyah
\cite{At} proves the additivity of the Adams operations defined in
terms of symmetric groups.

Now we set $\G =\{1\}$ and consider the diagonal embedding
$\Delta_n: X \rightarrow X^{\Delta} \hookrightarrow X^n$. Since
$\Hy_n$ acts on $X^{\Delta}$ trivially, it follows by
Proposition~\ref{prop_triv} that $\underline{K}^-_{\Hy_n}(X) \cong
\underline{K}(X) \otimes R^-(\Hy_n).$ We have the following
commutative diagram ($n$ odd):

\begin{eqnarray*}
  \begin{array}{rllll}
  \underline{K}(X)  &   & & &\\
   \boxtimes n \circ \aleph \downarrow \quad
  & \searrow \otimes n \circ \aleph
  &   && \\
 \underline{K}^-_{\Hy_n}(X^n)
 &\stackrel{\Delta_n^*}{\longrightarrow}
 & \underline{K}^-_{\Hy_n}(X) &\stackrel{\cong}{\longrightarrow}
   & \underline{K}(X) \otimes R^-(\Hy_n) \\
   {pr}\circ \phi_n \downarrow  \uparrow \iota
   &  & \iota
   \uparrow  \downarrow  {pr}\circ \phi^{\Delta}_n
  & \swarrow{ev} &\\
       \underline{K} \left( X^{\Delta} \right)
       &\stackrel{\vartheta}{\longrightarrow} &
      \underline{K} \left( X \right)
      & &
   \end{array}
\end{eqnarray*}
where $\phi^{\Delta}_n$ is the analog of $\phi_n$ when $X^n$ is
replaced by the diagonal $X$, and the evaluation map $ev$ is
defined to be the character value at the conjugacy class of type
$((n), \emptyset)$. The map from $\underline{K}(X)$ to itself
obtained along the left-bottom route in the above diagram
coincides with $\frac1n \psi^n$ given in
Definition~\ref{def_operat}. The map from $\underline{K}(X)$ to
itself obtained along the top-right route in the above diagram
gives $\frac1n$ times the $n$-th Adams operation. This of course
gives another proof that both $\psi^n$ and $\varphi^n$ are
additive when $\G$ is trivial.
\end{demo}

\begin{remark} \rm
$\tFG$ is a $Q$-$\lambda$ ring with $\varphi^n$ ($n$ odd) as the
$n$th Adams operation. If $X$ is a point and thus ${\cal
F}^-_{\G}(pt) = R^-_{\G},$ then our result reduces to the fact
that ${\cal F}^-_{\G}(pt)$ is a free $Q$-$\lambda$ ring generated
by $\G_*$. In particular when $\G$ is trivial this is isomorphic
to the model $Q$-$\la$ ring $\Omega_\C$.
\end{remark}

Denote by $\widehat{\cal F}_\G^t (X)$ the completion of $\tFGvar$
which allows formal infinite sums. Given $V \in
\underline{K}_\G(X)$, we introduce $Q(V, t) \in \widehat{\cal
F}_\G^t (X)$ as follows:
\begin{eqnarray} \label{eq_vo}
 Q(V, t) = \bigoplus_{n \geq 0} t^n (V|V)^{\boxtimes n}.
 \end{eqnarray}

The following lemma is immediate by Definition~\ref{def_operat}
and Remark~\ref{rem_tensor}.
\begin{lemma}
Given $V \in \underline{K}_\G(X)$, we have (for $r$ odd)

$$\varphi^r \circ \aleph (V) =\sum_{c \in \G_*} \zeta_c^{ -1}
          \phi_r^{-1}(\phi_r)_{c_r} (V^{\boxtimes r} \star \sigma_r
          (c)). $$
\end{lemma}

\begin{proposition} \label{prop_expon}
Given $V \in \underline{K}_\G(X)$, we can express $Q(V, t)$ as
follows:
\begin{eqnarray*}
 Q(V, t) & =& \exp \left( \sum_{r >0 \; odd}
 \frac2r \varphi^r \circ \aleph (V) t^r \right).
   \nonumber
 \end{eqnarray*}
Here the right-hand side is understood in terms of the algebra
structure on $\tFG$.
\end{proposition}

\begin{demo}{Proof}
By (\ref{eq_basicspin}) and the above lemma, we have

\begin{eqnarray*}
 Q(V, t) & =& \bigoplus_{n \geq 0} t^n (V^{\boxtimes n} \star \xi^n ) \\
   & =& \bigoplus_{n \geq 0} t^n
     \left( \sum_{\rho\in{\cal OP}_n(\G_*)} \phi_n^{-1}(\phi_n)_{\rho}
     (V^{\boxtimes n} \star 2^{l(\rho)} Z_{\rho}^{ -1}
     \sigma^{\rho})
     \right)                    \\
   & =& \bigoplus_{n \geq 0} \sum_{\rho\in{\cal OP}_n(\G_*)}
      \phi_n^{-1}(\phi_n)_{\rho}
      ( 2^{l(\rho)} Z_{\rho}^{ -1}
      t^n V^{\boxtimes n} \star \sigma^{\rho}
      )                   \\
   & =& \prod_{c \in G_*, r \geq 1\; odd}\frac1{m_r (c)!}
      \left( \frac2r t^r \zeta_c^{ -1}
      \phi_r^{-1}(\phi_r)_{c_r} (V^{\boxtimes r} \star \sigma_r (c))
      \right)^{m_r (c)}                                    \\
   & =&  \exp \left( \sum_{ r \geq 1\; odd} \frac 2r t^r
                   \sum_{c \in \G_*} \zeta_c^{ -1}
          \phi_r^{-1}(\phi_r)_{c_r} (V^{\boxtimes r} \star \sigma_r (c))
           \right)
             \\
   & =&  \exp \left(  \sum_{ r \geq 1\; odd}
   \frac 2r t^r \varphi^{r} \circ \aleph (V)
              \right) .
\end{eqnarray*}
\end{demo}

Combining with the additivity of $\varphi^r$, the proposition
implies
\begin{corollary}
 The following equations hold for $V, W \in \underline{K}_G(X)$:
 \begin{eqnarray*}
  Q( V \bigoplus W, t) & =& Q( V, t)Q( W, t) \\
     Q( -V, t) & =& Q(V, -t).
 \end{eqnarray*}
\end{corollary}

\begin{remark}  \rm
The generating function $Q(V, t)$ is essentially half the twisted
vertex operator, and the other half can be obtained by the adjoint
operator to $Q(V, t)$. Twisted vertex operators have played an
important role in the representation theory of infinite
dimensional Lie algebras and the moonshine module, cf. \cite{FLM}.
When $X$ is a point, we can develop the picture more completely
(cf. \cite{JW}) to provide a group theoretic realization of vertex
representations of twisted affine and twisted toroidal Lie
algebras (also compare \cite{FJW2} for a different construction).
\end{remark}
\subsection{Twisted Heisenberg algebra and $\tFG$} \label{sec_heis}

We see from Theorem~\ref{th_size} that $\tFG$ has the same size of
the tensor product of the Fock space of an infinite-dimensional
twisted Heisenberg algebra of rank $\dim K^0_\G(X)$ and that of an
infinite-dimensional twisted Clifford algebra of rank $\dim
K^1_G(X)$. In this section we will actually construct such a
Heisenberg/Clifford algebra, which we will simply refer to as a
twisted Heisenberg (super)algebra from now on.

The dual of $\underline{K}_\G(X)$, denoted by
$\underline{K}_\G(X)^*$, is naturally $\Z_2$-graded as identified
with $\underline{K}^0_\G(X)^* \bigoplus \underline{K}^1_\G(X)^*$.
Denote by $\langle \cdot, \cdot \rangle$ the pairing between
$\underline{K}_\G(X)^*$ and $\underline{K}_\G(X)$. For any $m \in
2\Z_+ +1 $ and $\eta \in \underline{K}_\G(X)^*$, we define an
additive map
\begin{eqnarray}  \label{eq_ann}
  a_{-m} (\eta) : \underline{K}_{G_n}(X^n) \longrightarrow
            \underline{K}_{G_{n-m} }(X^{n-m} )
\end{eqnarray}
as the composition
\begin{eqnarray*}
  \underline{K}^-_{\tG_n}(X^n)
& \stackrel{Res}{\longrightarrow} & \underline{K}^-_{\tG_m
\widehat{\times} \tG_{n -m} }(X^n)
\stackrel{k^{-1} }{\longrightarrow}
              \underline{K}^-_{\tG_{m}} (X^{m} )\bigotimes
               \underline{K}^-_{\tG_{n -m} }(X^{n -m} )     \\
&\stackrel{{ch}_m \otimes 1 }{\longrightarrow} &
\underline{K}_{\G} (X ) \bigotimes \underline{K}_{\tG_{n -m}}
(X^{n -m} )
\stackrel{\eta \otimes 1 }{\longrightarrow}
              \underline{K}^-_{\tG_{n -m}} (X^{n -m} ).
\end{eqnarray*}

On the other hand, we define for any $m \in 2\Z_+ +1$ and $V \in
\underline{K}_\G (X)$ an additive map
\begin{eqnarray}  \label{eq_creat}
  a_{m} (V) : \underline{K}^-_{\tG_{n-m} }(X^{n-m} )
          \longrightarrow \underline{K}^-_{\tG_n}(X^n)
\end{eqnarray}
as the composition
\begin{eqnarray*}
  \underline{K}^-_{\tG_{n-m} }(X^{n-m} )
  &\stackrel{\frac{m}2\varpi_m (V) \boxtimes \cdot}{\longrightarrow} &
  \underline{K}^-_{\tG_m }(X^m)
  \bigotimes \underline{K}^-_{\tG_{n -m}} (X^{n -m} ) \\
  & \stackrel{k }{\longrightarrow} &
    \underline{K}^-_{\tG_m \widehat{\times} \tG_{n -m}}(X^n )
\stackrel{\mbox{Ind}}{\longrightarrow}
          \underline{K}^-_{\tG_n}(X^n).
\end{eqnarray*}

Let $\cal H$ be the linear span of the operators $a_{-m} (\eta),
a_{m} (V), m \in 2\Z_+ +1$, $\eta \in \underline{K}_\G(X)^*,$ $V
\in \underline{K}_\G(X)$. Clearly $\cal H$ admits a natural
$\Z_2$-gradation induced from that on $\underline{K}_\G(X)$ and $
\underline{K}_\G(X)^* $. Below we shall use $[ - , - ]$ to denote
the supercommutator as well. It is understood that $[a, b]$ is the
anti-commutator $ab +ba$ when $a, b \in \cal H$ are both odd
elements according to the $\Z_2$-gradation.

\begin{theorem}    \label{th_heisenberg}
When acting on $\tFG$, $\cal H$ satisfies the twisted Heisenberg
superalgebra commutation relations, namely for $m , l \in
2\Z_++1,$ $ \eta, \eta ' \in \underline{K}_\G(X)^*, $ $V,W \in
\underline{K}_\G(X)$, we have
\begin{eqnarray*}
  [ a_{-m} (\eta), a_{l} (V)]
  & =& \frac l2 \delta_{m,l} \langle \eta, V \rangle ,
           \label{eq_heis}  \\
  {[ a_{m} (W)    , a_{l} (V)]} & =& 0,      \label{eq_induct} \\
  {[ a_{-m} ({\eta}), a_{-l} ({\eta}') ] }& = & 0.  \label{eq_restrict}
\end{eqnarray*}
Furthermore, $\tFG$ is an irreducible representation of the
twisted Heisenberg superalgebra.
\end{theorem}

\begin{remark}  \rm
The proof of the Heisenberg algebra commutation relation can be
given in a parallel way as the one for Theorem 4, \cite{W1}. The
irreducibility of $\tFG$ as module over the Heisenberg algebra
follows from Theorem~\ref{th_size}. Given a bilinear form on
$\underline{K}(X)$, then we can get rid of $\underline{K}(X)^*$ in
the formulation of the above theorem. In the special case when $X$
is a point, the Heisenberg algebra here specializes to the one
given in \cite{JW} acting on $R^-_{\G}$. In the case when $\G$ is
a finite subgroup of $SL_2 (\C)$, we may consider further the
space which is the tensor product of ${\cal F}^-_{\G}(pt)$ with a
module of a certain $2$-group which can be constructed out of
$\G$, and realize in this way a vertex representation of a twisted
affine and a twisted toroidal Lie algebra. This is treated in
\cite{JW} in detail.
\end{remark}
\section{Appendix: Another formulation using $\widetilde{S}_n$ and
$\widetilde{\G}_n$}

As is well known (cf. e.g. \cite{Sc, J, Jo, HH}), the symmetric
group $S_n$ has a double cover $\widetilde{S}_n$:

\begin{eqnarray*}
1 \longrightarrow\Z_2 \longrightarrow \widetilde{S}_n
\stackrel{\theta_n}{\longrightarrow} S_n \longrightarrow 1,
\end{eqnarray*}
generated by $z$ and $ t_i, i=1, \cdots, n-1$ and subject to the
relations: $$
 z^2=1,  \quad t_i^2=(t_it_{i+1})^3=z, \quad
 t_it_j=zt_jt_i \;(i>j+1), \quad zt_i=t_iz.
$$
The map $\theta_n$ sends $t_i$'s to the simple reflections
$s_i$'s in $S_n$. The group $\widetilde{S}_n$ carries a natural
$\Z_2$ grading by letting $t_i$'s be odd and $z$ be even.

Given a finite group $\G$, the symmetric group $S_n$ acts on the
product group $\G^n$, and the group $\widetilde{S}_n$ acts on
$\G^n$ via $\theta_n$. Thus we can form a semi-direct product
$\widetilde{\G}_n := \G^n \rtimes_{\theta_n} \widetilde{S}_n$,
which carries a natural finite supergroup structure by letting
elements in $\G^n$ be even, cf. \cite{FJW2}. We still denote by
$\theta_n$ the quotient map $\widetilde{\G}_n \rightarrow
\widetilde{\G}_n/\langle 1,z \rangle =\Gn$.

Given a $\G$-space $X$, we have seen $X^n$ affords a natural $\Gn$
action. Then we can apply the general construction in
Sect.~\ref{subsect_ktwist} to construct the category ${\cal
C}^-_{\widetilde{\G}_n}(X^n)$ and its associated $K$-group
$K^-_{\widetilde{\G}_n}(X^n)$. As before, we denote
$\underline{K}^-_{\widetilde{\G}_n}(X^n)
=K^-_{\widetilde{\G}_n}(X^n) \otimes \C$.

We then form the direct sum
\begin{eqnarray*}
 {\frak F}^-_\G(X) &=&
 \bigoplus_{n=0}^{\infty} \underline{K}^-_{\widetilde{\G}_n}(X^n)
 \\
 {\frak F}^-_\G(X, t) &=&
 \bigoplus_{n=0}^{\infty} t^n \underline{K}^-_{\widetilde{\G}_n}(X^n).
\end{eqnarray*}
When $X$ is a point, $\underline{K}^-_{\widetilde{\G}_n}(pt) $
reduces to the Grothendieck group $R^-(\widetilde{\G}_n)$ of spin
supermodules of $\widetilde{S}_n$, and ${\frak F}^-_\G(pt)$ has
been studied in detail in \cite{FJW2}. The purpose of the Appendix
is to outline how to modify the various constructions of algebraic
structures on $\tFG$ for the new space ${\frak F}^-_\G(X)$. As the
constructions are very similar to the $\tFG$ case, we will be
rather sketchy.

Given $n, m \ge 0$, we can define a $\Z_2$-graded subgroup
$\widetilde{\G}_n \widehat{\times} \widetilde{\G}_m$ of
$\widetilde{\G}_{n+m}$, in a way analogous to (\ref{eq_product}).
Then the obvious analog of constructions (\ref{eq_mult}) and
(\ref{eq_comult}) defines a multiplication and comultiplication on
the space ${\frak F}^-_\G(X)$. The following is an analog of
Theorem~\ref{th_hopf} and it generalizes Theorem 3.8 of
\cite{FJW2} which is our special case when $X$ is a point.

\begin{theorem}
The space  ${\frak F}^-_\G(X)$ carries a natural Hopf algebra
structure.
\end{theorem}

By the analysis of the split conjugacy classes given in the proof
of Theorem 2.5 of \cite{FJW2}, we see that the
$\widetilde{\G}_n$-space $X^n$ satisfies the strong vanishing
property, and the analog of Lemma~\ref{lem_orbit} holds. Therefore
we obtain the following theorem which is an analog of
Theorem~\ref{th_size}.

\begin{theorem} \label{th_newsize}
As a $(\Z_+ \times \Z_2)$-graded algebra ${\frak F}^-_{\G}(X, t)$
is isomorphic to the supersymmetric algebra
  $ {\cal S} \left( \bigoplus_{ r=1}^{\infty} t^{2r-1} \underline{K}_\G(X)
      \right)$.
\end{theorem}

Except the first two terms and the first two arrows in the diagram
(\ref{eq_maps}), the rest of the diagram has a direct analog for
$K^-_{\widetilde{\G}_n}(X^n)$. Note in the definition of the
$K$-theory maps $ch_n$ and $\varpi_n$ (see
Definition~\ref{def_operat}) only the part of the diagram
(\ref{eq_maps}) which can be directly generalized to the
$K^-_{\widetilde{\G}_n}(X^n)$ setup has been used. Therefore
analog of $ch_n$ and $\varpi_n$ can be defined in our new setup.
This guarantees the analog of annihilation operators
(\ref{eq_ann}) and the creation operators (\ref{eq_creat}) can be
defined in our new setup. In this way we obtain the following
which is an analog of Theorem~\ref{th_heisenberg}.

\begin{theorem}
The space  ${\frak F}^-_\G(X)$ affords an action of the twisted
Heisenberg algebra $\cal H$ in terms of natural additive
$K$-theory maps. Furthermore this representation is irreducible.
\end{theorem}

We remark that it is much less natural to use $\widetilde{S}_n$ to
construct various $K$-theory operations on $K(X)$ as done in
Sect.~\ref{sec_oper} using a double cover $\widetilde{H}_n$ of the
hyperoctahedral group.

The connection between $\tFG$ and the $Q$-$\la$ ring in
Sect.~\ref{subsect_twistvo} does carry over to our new setup.
Keeping Remark~\ref{rem_tensor} in mind and knowing that
$\widetilde{S}_n$ also has a so-called basic spin supermodule (cf.
e.g. \cite{Jo, FJW2}), we can use it to define the analog of
(\ref{eq_vo}). Indeed this also defines an analog of the map
$\boxtimes^n \circ \aleph$ (cf. the diagram (\ref{eq_maps})), and
thus an analog of $\varphi^n \circ \aleph$. Therefore, we have an
analog of Proposition~\ref{prop_expon} in our new setup which
generalizes Proposition~6.2 in \cite{FJW2}. However there is an
unpleasant square root of $2$ in the formula which originates in
the spin representation theory of $\widetilde{S}_n$ and
$\widetilde{\G}_n$. This is another reason why we have preferred
the formulation in the main body of the paper using
$\widetilde{H}_n$ and $\tG_n$.

One may wonder that why $K^-_{\tG_n}(X^n)$ and
$K^-_{\widetilde{\G}_n}(X^n)$ are so similar to each other and
there are almost parallel constructions on $ \tFG$ and ${\frak
F}^-_\G(X)$. When $X$ is a point and thus the $K$-groups reduces
to the corresponding Grothendieck groups of spin supermodules,
this has been noticed by various different authors (cf. e.g.
\cite{Ser, Jo2, Naz, Y, FJW2} and the references therein).
Yamaguchi \cite{Y} explains clearly such a phenomenon by
establishing an isomorphism between the group superalgebra $\C
[\widetilde{H}_n]/ \langle z=-1 \rangle$ and the ($\Z_2$-graded)
tensor product of the group superalgebra $\C[\widetilde{S}_n]/
\langle z=-1 \rangle$ with the complex Clifford algebra $C_n$ of
$n$ variables (this is not the same copy of $C_n$ associated to
the subgroup $\NP_n$ in $\widetilde{H}_n$!). It follows that the
group superalgebra $\C [\tG_n]/ \langle z=-1 \rangle$ is
isomorphic to the tensor product of the group superalgebra
$\C[\widetilde{\G}_n]/ \langle z=-1 \rangle$ with $C_n$. Note that
a Clifford algebra admits a unique irreducible supermodule. As
this $C_n$ acts on an $\tG_n$-bundle over $X^n$ fiberwise, this
isomorphism provides a direct isomorphism between
$K^-_{\tG_n}(X^n)$ and $K^-_{\widetilde{\G}_n}(X^n)$.

We can also forget about the $\Z_2$-gradings (i.e. the super
structures) in the group $\widetilde{\G}_n$, in the construction
of the category ${\cal C}^-_{\widetilde{\G}_n}(X^n)$ and its
associated $K$-group $K^-_{\widetilde{\G}_n}(X^n)$. Let us denote
the resulting new $K$-group by $K^s_{\widetilde{\G}_n}(X^n)$. In
particular when $X$ is a point and $\G$ is trivial, this reduces
to the Grothendick group $R^s(\widetilde{\G}_n)$ of spin (not
super) modules of $\widetilde{S}_n$ where $z$ still acts as $-1$.
We can then apply the decomposition theorem of Adem-Ruan \cite{AR}
to calculate $K^s_{\widetilde{\G}_n}(X^n) \otimes \C$ in terms of
$K_\G(X) \otimes \C$. The difference here from the calculations in
Theorem~\ref{th_size} and Theorem~\ref{th_newsize} is that the
{\em odd} split conjugacy classes of $\widetilde{\G}_n$ will also
make contributions. Recall that the orbifold Euler number
$e(X,\G)$ defined in \cite{DHVW} is the same as the Euler number
of the equivariant $K$-theory $K_\G(X)$. Using the description of
{\em even} and {\em odd} split conjugacy classes of
$\widetilde{\G}_n$ (cf. \cite{FJW2}, Theorem~2.5), we can obtain
the Euler number of $K^s_{\widetilde{\G}_n}(X^n)$, denoted by
$e^s(X^n, \widetilde{\G}_n)$, in terms of the following generating
function (compare Corollary~\ref{cor_euler}):

\begin{eqnarray*}
 \sum_{n=0}^{\infty} t^n  e^s(X^n,\widetilde{\G}_n)
 &=&  { \prod_{ r=1}^{\infty}
    (1 - t^{2r-1})^{-e(X,\G)} } \\
 &&+ \prod_{ r=1}^{\infty}
 (1 + t^{2r-1})^{e(X,\G)} \cdot \\
 &&\quad \cdot \frac12 \ \left( \prod_{ r=1}^{\infty}
 (1 + t^{2r})^{e(X,\G)} -\prod_{ r=1}^{\infty}
 (1 - t^{2r})^{e(X,\G)} \right).
\end{eqnarray*}

The second summand in the right-hand side of the above equation
counts the contributions from odd split conjugacy classes. When we
set $X$ to be a point and $\G$ trivial (and thus $e(X,\G) =1$),
this formula reduces to the classical generating function for the
spin Grothendick group $R^s(\widetilde{S}_n)$ (cf. Theorem~3.6,
\cite{Jo}, pp.~213; Corollary~3.10, \cite{HH}, pp.~32). We remark
that due to some inaccurate analysis of split conjugacy classes of
$\widetilde{S}_n$, the formula (6.10) given in \cite{Di} (even in
the case when $\G$ is trivial and $X$ is a point) is incompatible
with this classical statement.

Department of Math., North Carolina State University, Raleigh, NC
27695.

{\em Current address}: Department of Mathematics, University of
Virginia, Charlottesville, VA 22904, U.S.A. Email:
ww9c@virginia.edu

\end{document}